\documentclass[12pt]{article} 
\usepackage{graphicx} 
\usepackage{graphics}
\usepackage{amsfonts} 
 \usepackage{amsthm,amssymb} 
\usepackage{latexsym,amsfonts,amssymb,amsmath,euscript} 
\def\E{I \kern -0.37 em E}      
 
\def\pr{\mathbb{P}} 
\def\N{I \kern -0.37 em N}      %
\def\D{I \kern -0.37 em D}      %
\def\R{I \kern -0.37 em R}      %
\def\P{I \kern -0.37 em P}

\def\v{\bigvee}

\def\bm{\vspace{05mm} \begin{displaymath}}  
\def\be{\vspace{05mm} \begin{equation}}  
\def\ee{\end{equation} \vspace{05mm}}  
\def\v5{\vspace{05mm}}

\def\fedum{i}\def\feddois{ii}\def\fedtres{iii} 
\def\assone{1}\def\asstwo{2}\def\assthree{3}\def\assfour{4}
\reversemarginpar 
\newcounter{mmm}  

\newcounter{rnb}  
\newcommand{\remlabel}[1]{\refstepcounter{rnb}\label{#1}\thernb} 
\newcounter{fnb}  

\newcounter{lem}  
\newcommand{\corlabel}[1]{\refstepcounter{lem}\label{#1}\thelem} 
\newcounter{htq}  

\newcounter{trw}  

\newcounter{fmq}  

\newcounter{fqx}  

\newcounter{mish}  
\newcommand{\mishlabel}[1]{\refstepcounter{mish}\label{#1}\themish} 
\newcounter{ammel}  
\newcommand{\ammellabel}[1]{\refstepcounter{ammel}\label{#1}\theammel} 
\newcounter{awk}   

\newcounter{conj}   
\newcommand{\conjlabel}[1]{\refstepcounter{conj}\label{#1}\theconj} 
\def\price{p}  
\def\demand{d}  
\def\tprice{\tilde\price}
\def\tdemand{\tilde\demand} 
\def\feedback{\lambda} 
\def\medaval{\bar{v}} 
\def\ve{v} 
\def\bard{\bar{d}} 
\def\proportion{\alpha} 
\def\decide{d} 
\def\sus{J} 
\def\ini{{\rm initial}}  
\def\fini{{\rm final}} 
\def\diftza{a}
\def\maxdir{b}
\def\mondir{c}
\def\eqibda{d}
\def\expez{e}
\def\finvar{f}
\begin{document}      
\begin{center} 
{\LARGE \bf Stability analysis with applications of a two-dimensional 
dynamical system arising from a stochastic model for an asset market} 
\\ 
\mbox{\kern1em}\\ 
 Vladimir Belitsky, Antonio Luiz Pereira\\  
{\em Institute of Mathematics and Statistics, University of S\~ao Paulo, Brazil} \\ 
Fernando Pigeard de Almeida Prado\\ 
{\em Departamento de F\'\i sica e Matem\'atica, FFCLRP, Universidade de S\~ao
Paulo, Brazil} 
\end{center}  
 
\begin{center} {\bf Abstract}\end{center}{\small  
We analyze the stability properties of equilibrium solutions and periodicity of 
orbits in a two-dimensional dynamical system whose orbits mimic the evolution 
of the price of an asset and the  excess demand for that asset. The  construction 
of the system  is grounded upon a heterogeneous interacting agent model for a  
single risky asset market. An advantage of this construction procedure is that 
the resulting dynamical system becomes a macroscopic market model which mirrors 
the market quantities and qualities that would typically be taken into account 
solely at the microscopic level of modeling. The system's parameters correspond 
to: (a) the proportion of speculators in a market; (b) the traders' speculative 
trend; (c) the degree of heterogeneity of idiosyncratic 
evaluations of the market agents with  respect to the asset's fundamental value; 
and (d) the strength of the feedback of the population excess demand on the asset price 
update increment. This correspondence allows us to employ our results in order 
to infer plausible causes for the emergence of  price and demand fluctuations 
in a  real asset market.   
 
The employment of dynamical systems for studying evolution of stochastic models 
of socio-economic phenomena is quite usual in the area of heterogeneous interacting 
agent models. However, in the vast majority of the cases present in the literature, 
these dynamical systems  are one-dimensional. Our work is among the few in the area 
that construct and study two-dimensional dynamical systems and apply them for 
explanation of socio-economic phenomena.
} 
 
 
\medskip\noindent{\sl Classifications:} 
\newline AMS classification numbers (MSC2000):  
60J20, 
60K35,  
82C22, 
58F08. 
\newpage

\section{Introduction}\label{introduction} 
We investigate equilibria and stability properties of the discrete time dynamical system 
generated by  planar map $\Psi$ in the following manner 
\begin{equation}\label{DYMA}\left( \begin{array}{c} 
               \price_{n} \\ 
                \demand_{n} 
                 \end{array} \right)= 
\Psi\left( \begin{array}{c} 
               \price_{n-1} \\ 
                \demand_{n-1} 
                 \end{array} \right) =  
                   \left( \begin{array}{l} 
               \price_{n-1}+\feedback\demand_{n-1} \\ 
                \proportion\big[1-2\Phi\big(\price_{n-1}+(\feedback-\sus)\demand_{n-1}\big)\big]\\ \mbox{\kern5em}+(1-\proportion) 
\big[1-2\Phi\big(\price_{n-1}+\feedback\demand_{n-1}\big)\big]                 \end{array} \right). 
\end{equation} 
In the dynamical system (\ref{DYMA}), $\proportion\in (0,1)$, $\feedback$ and $\sus$ are 
positive real numbers and $\Phi$ is a probability distribution function. Since we do not 
impose rigid constraints on this distribution function, the analysis of (\ref{DYMA}) 
becomes a nontrivial task. Nevertheless, our  study of (\ref{DYMA}) was motivated by 
its potential applications rather than by intrinsic mathematical challenges yielded by the 
non-linearity. The applications stem from the link of the system to a stylized  model of 
a single risky asset market represented by a stochastic process constructed by us.
 
The stochastic process just mentioned is a microscopic model of a market in the sense that
it mirrors individual behavior of each market agent. The model will be explained in 
Section~\ref{model}. We constructed it 
using  ideas from the area of Heterogeneous Interacting Agent Models\footnote{The works  
          \cite{glaeser-scheinkman} and \cite{judd-tesfatsion} provide an broad 
          survey of the history of the HIAMs, as well as of the current   
          state-of-the-art in the area. In particular, our HIAM  shares many  features in common 
          with the HIAM from the seminal works of Kirman (\cite{kirman-ants}) and of Lux (\cite{lux}). 
        Among the more recent studies that present HIAMs similar to ours, we wish to single out  \cite{gordon-seller} and \cite{nadal-multiple-equilibria} because the problems addressed by them 
        are related to the issues discussed by us here.  
          However, we shall not pursue a detailed comparison of our model to the already existing ones  
          since this  is not essential for our presentation.} (HIAMs).    
HIAMs have a long history, wide applications and employ methods from a variety
of  disciplines including  Statistical Mechanics and Interacting Particle Systems 
(see \cite{liggett}). The experience shows\footnote{See the review works  
\cite{glaeser-scheinkman} and \cite{judd-tesfatsion}, and citations  therein.}   
that HIAMs are able to  adequately describe macroscopic characteristics of 
socio-economic systems by means of mirroring system's microscopic components and 
their mutual interactions. Our HIAM was constructed with the aim 
to investigate evolution of two particular macroscopic characteristics: 
the asset price and the  population excess demand for an asset in real world 
markets. In the constructed model, the evolutions of interest are represented 
by stochastic processes. These processes possess the following property: as 
the model's agent population increases, their trajectories converge to orbits of the 
dynamical system (\ref{DYMA}) (this convergence will be explained in 
Section~\ref{dynas}). This convergence yields the following interpretation of the parameters 
and variables of (\ref{DYMA}) in terms of a single risky asset market: 
\begin{center}
\begin{tabular}{cl}
$\price_n$ & corresponds to the price of the market asset at time $n$;\\
$\demand_n$ & corresponds to the population excess demand for the asset at time $n$\\
            & (in the sequel, it will be called excess demand for short);\\
$\proportion$ & corresponds to the proportion of speculators among the market traders;\\
$\sus$ & corresponds to the traders' speculative trend;\\
$\feedback$ & corresponds to the feedback of the excess demand on the increment\\
            &  of the update of the asset price;
\end{tabular}\end{center}
\begin{center}
\begin{tabular}{cl}
$\Phi$ & corresponds to the distribution of the deviations around a constant $\medaval$\\ 
       &  of individual evaluations of the asset fundamental value, where $\medaval$ corresponds\\ 
       & to the market fundamental value of the asset.
\end{tabular}\end{center}
This interpretation will be explained in details and justified in 
Sections~\ref{model} and \ref{dynas}. 
It will allow us to re-phrase our results so that they explain  whether and why an asset  
price and excess demand in an asset market would or would not oscillate as time goes on; this will 
be the contents of Section~\ref{application}.
 
The above mentioned  convergence of HIAM's trajectories to orbits of an appropriate 
dynamical system is a well known fact in the HIAM area. It is usually 
employed,\footnote{See  \cite{levy-social} for a general account on the 
         employment of this convergence for revealing socio-economic phenomena.}
but almost always in such a way that  the  resulting dynamical
system is of the dimension one. This is motivated by the wish to 
use the graphical analysis as a tool for study this system. 
Contrasting, in the present work, the 
limiting dynamical system is two-dimensional. The dimension increase happens to make 
a difference because it allows one to see socio-economic phenomena that are invisible through 
the lenses of one-dimensional dynamical systems. The price to pay for this is the growing
complexity of the proofs.
 
The results of our study of (\ref{DYMA}) are formulated in Theorems~\ref{propdyns}, \ref{limo}  
and \ref{dyns} of  Section~\ref{results}. Theorem~\ref{propdyns} asserts  that the 
point $(0,0)$ is the unique  equilibrium state of the dynamical system (\ref{DYMA}) 
and determines when  it is locally asymptotically stable and when it is
unstable.  An interesting feature of this theorem is that it characterizes the 
stability/instability via a relation between just two expressions involving the 
model's parameters (that are, we recall, three real numbers and a probability 
distribution function).  This reduction of the parameters' space allows us to 
draw a phase diagram (Figure~\ref{regions}) representing the stability/instability of 
the unique equilibrium point.  Our second result, Theorem~\ref{limo}, gives conditions  
on the parameter values for the appearance of a stable periodic orbit of the dynamical 
system. Generalizations of this theorem are indicated in Conjectures~\ref{noche} and 
\ref{checked}. Our third result, Theorem~\ref{dyns}, identifies  parameter values for which 
the point $(0,0)$ is a globally stable equilibrium. 

\section{Results}\label{results}
Here, we present results of our study of the dynamical system (\ref{DYMA}) in which
\begin{equation}\label{parcon}
\proportion\in (0,1),\,\,\sus>0,\,\, \feedback>0,
\end{equation}
and $\Phi$ is a probability distribution function that satisfies the following conditions:
\begin{equation}\label{phicon}\begin{array}{cl}
(\diftza) & 
\Phi \hbox{ is everywhere differentiable; we shall denote its derivative by }\Phi^\prime;\\
(\maxdir) & 
\sup_x\left\{\Phi^\prime(x)\right\}=\Phi^\prime(0);\\
(\mondir) & 
\Phi^\prime\hbox{ is strictly increasing on }(-\infty, 0),\hbox{ and strictly decreasing on }(0,\infty);\\
(\eqibda) & 
\Phi(0)=1/2;\\
(\expez) & 
\hbox{a random variable with the distribution }\Phi\hbox{ has zero expectation};\\
(\finvar) & 
\hbox{a random variable with the distribution }\Phi\hbox{ has finite variance.}\end{array}
\end{equation}
We might have relaxed significantly the constraints (\ref{phicon}) on $\Phi$, but this 
would not have  broadened  the  application of our results and  would make their proofs 
more cumbersome  without introducing  essentially new ideas. We note also that for 
applications, it is reasonable to admit that $\Phi$ is a Gaussian zero mean distribution function.
This function satisfies the above conditions.

\subsection{The uniqueness of equilibrium and its basic stability properties}\label{osnov}

\noindent{\bf Theorem~\mishlabel{propdyns}} \ (uniqueness of equilibrium and its stability). 
{\em  
\begin{description}
\item{\bf (a)}  The origin (i.e., the point $(0,0)$) is the unique equilibrium of the dynamical system
(\ref{DYMA}).
\item{\bf (b)} If  $\feedback  \Phi^\prime(0) > 1+ 2 \proportion \sus \Phi^\prime(0))$,
then   the origin is an unstable equilibrium of (\ref{DYMA}), while if
$\feedback \Phi^\prime(0) \leq 1+ 2 \proportion \sus \Phi^\prime(0)$, then the stability
 of the origin depends on the value of $2 \proportion \sus\Phi^\prime(0)$ in the following manner:
\begin{description}\item{} if $ 2  \proportion \sus\Phi^\prime(0)  <  1$
then the origin is locally     asymptotically stable;
   \item{} if $ 2  \proportion \sus\Phi^\prime(0) > 1$ then the origin is unstable.
\end{description}
\end{description}
}
  
\bigskip 
\noindent {\em Proof of Thm~\ref{propdyns}.} \  If $ (\price_n, \demand_n) = (\price_0,\demand_0)$, 
for $n \in \mathbb{N}$, then in virtue of the first equation of (\ref{DYMA}), $\feedback\demand_0=0$. 
This implies that $\demand_0 =0$ because $\feedback >0$ by construction. On substituting $\demand_0 =0$
in the second equation  of (\ref{DYMA}), we get that $1-2\Phi \big(\price_0 \big)  = 0$. 
This can  hold only if $\price_0=0$ because $\Phi$ is monotone and $\Phi(0)=1/2$ (the properties 
ensured by the assumptions (\ref{phicon})). This completes the proof of item (a). 
     
\begin{figure}[h] 
 \begin{center}
\begin{picture}(400,400)
\put(50,50){\vector(1,0){300}} 
\put(50,50){\vector(0,1){300}} 
\put(50,150){\line(1,1){200}} 
\put(150,50){\line(0,1){200}} 
\qbezier(50,100)(50,150)(150,250) 
\qbezier(150,250)(190,290)(280,350) 
\qbezier(50,100)(50,50)(150,50) 
\qbezier(150,50)(200,50)(350,80) 
\multiput(50,50)(100,0){3}{\circle*{2}} 
\multiput(50,100)(0,50){5}{\circle*{2}}  
\put(30,30){0}\put(149,30){1}\put(249,30){2}\put(320,30){$u=2\alpha J\Phi^\prime(0)$}
\put(30,99){1}\put(30,149){2}\put(30,199){3}\put(30,249){4}\put(30,299){5}
\put(60,360){$w=2\lambda\Phi^\prime(0)$} 
\put(65,30){\oval(20,20)[bl]}\put(55,30){\vector(0,1){35}}\put(70,20){Region 2} 
\put(80,150){\oval(40,40)[bl]}\put(60,150){\vector(0,1){5}}\put(80,130){Region 3} 
\put(90,80){Region 1a} 
\put(160,80){Region 1b} 
\put(230,30){\oval(20,20)[br]}\put(240,30){\vector(0,1){27}}\put(180,20){Region 4} 
\put(100,260){Region~6} 
\put(200,280){Region 5}\put(230,290){\vector(-1,1){18}} 
\qbezier(90,175)(150,170)(200,250)\put(90,175){\circle*{4}}\put(200,250){\circle*{4}}
\put(75,160){$(u_{{\rm initial}}, w_{{\rm initial}})$}  
\put(200,230){$(u_{{\rm final}}, w_{{\rm final}})$} 
%
%
\rotatebox{10}{\put(270,20){$w=(1-\sqrt{u})^2$}}
\rotatebox{30}{\put(250,200){$w=(1+\sqrt{u})^2$}}
\rotatebox{45}{\put(100,330){$w=2u+2$}}
\put(-200,150){\rotatebox{90}{$u=1$}}
\end{picture}
\end{center}
\centering
\caption{Phase diagram representing the stability of the unique equilibrium of the dynamical system (\ref{DYMA}).}\label{regions}
\end{figure}
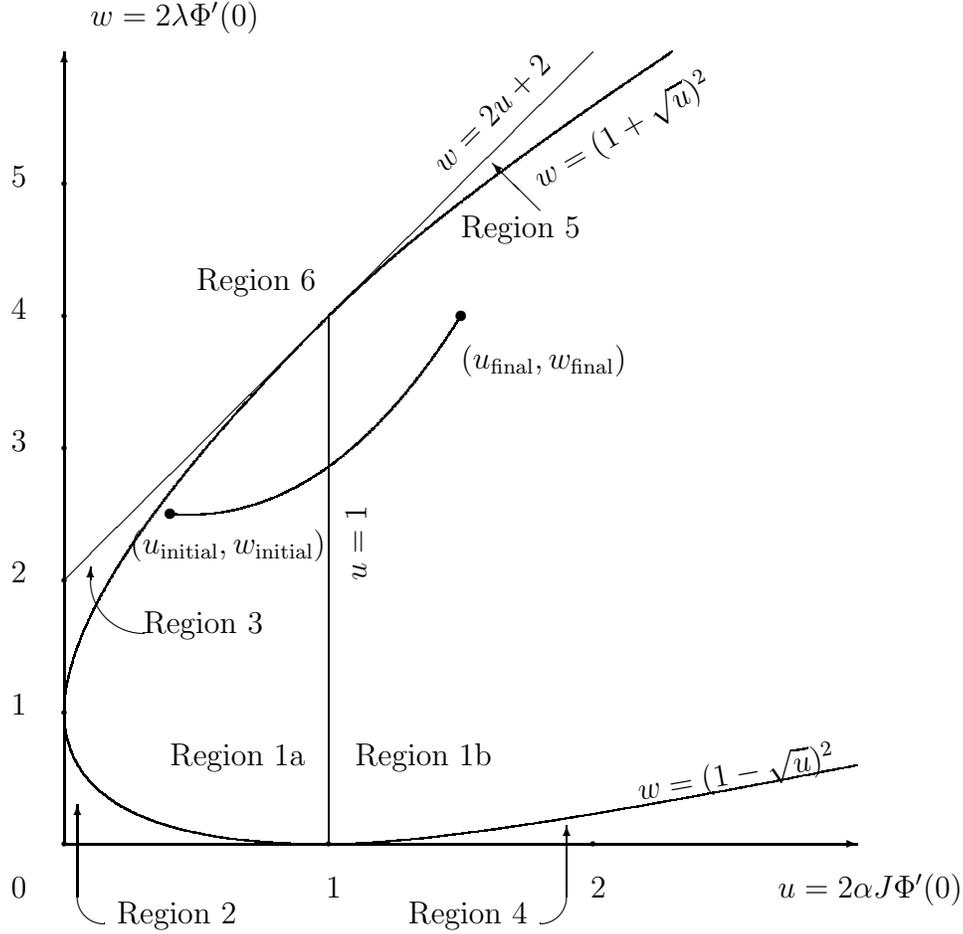

We proceed with the proof of (b) of the theorem. {}From (\ref{DYMA}) we easily get the Jacobian 
matrix of $\Psi$ at the origin: 
     \begin{equation} \label{jacob} 
     \left( \begin{array}{cc}
      1 & \feedback \\
      -2\Phi'(0) & 2 (\alpha\sus - \feedback) \Phi'(0)
      \end{array}               
     \right).
     \end{equation}
Its characteristic polynomial is
$
\mu^2 - \big[1 + 2 (\proportion \sus - \feedback) \Phi'(0)\big]\mu + 
2 \proportion\sus  \Phi'(0)$
 whose roots  are
\begin{equation}\label{rutik}\begin{array}{rcl}
        \mu_1 & = &  \frac{1}{2} \big[ 1 + 2 (\proportion \sus - \feedback) \Phi'(0) + \sqrt{\Delta} \big],\\ 
        \mu_2 & = &  \frac{1}{2} [ 1 + 2 (\proportion \sus - \feedback) \Phi'(0) - \sqrt{\Delta} \big],
\end{array}
\end{equation}      
where
\begin{equation}\label{Delta}
\Delta :=\left[1 + 2 (\proportion\sus - \feedback) \Phi'(0)\right]^2 - 8 \proportion\sus \Phi'(0).
\end{equation}
The assertions in (b) will all follow from  the principle of linearized stability 
via analysis of the values of $\mu_1$ and $\mu_2$.  The analysis are split into six cases. 
They may be seen from  Fig.~\ref{regions}. This figure is a kind
 two dimensional   phase diagram  because -- due to the form of the roots -- the stability analysis relies 
on  the relation between just two expressions: $2\proportion\sus  \Phi'(0)$ and $2 \feedback \Phi'(0)$.

All-throughout below, we shall usually use in our calculations that $\feedback>0$ (the constraint on $\feedback$ imposed by our construction) and that $\Phi^\prime(0)>0$ (the inequality that follows from (\ref{phicon}-\diftza, \maxdir, \mondir)). 
       
\medskip     
\noindent{\bf  Case 1:   $\Delta  \leq 0$} (region~1 in Fig.~\ref{regions}).
  
We have:
\begin{equation}\label{abis}
      |\mu_1|^2 = |\mu_2|^2  =   
      \frac{1}{4} (1 + 2 (\proportion\sus - \feedback) \Phi'(0))^2 + \left(\sqrt{-\Delta}\right)^2 
       =  2 \proportion\sus \Phi'(0), 
\end{equation}
and therefore the origin is locally asymptotically  stable, in case $2 \proportion\sus \Phi'(0) < 1$ 
(region~1a), and is unstable, in case $2 \proportion\sus \Phi'(0) > 1$ (region 1b).

\medskip
\noindent{\bf Case 2:   $\Delta > 0$, 
  $2 \feedback \Phi' (0) < 1 + 2 \proportion\sus \Phi'(0) $
        and $2 \proportion\sus \Phi'(0) < 1$} (region 2 in Fig.~\ref{regions}).
       
From the hypotheses, we obtain immediately:  $ 0 < 
          1+  2 (\proportion\sus -  \feedback) \Phi'(0)) <  1+  2  \proportion\sus  \Phi'(0)$.
By adding and subtracting $\big(1+2\proportion \sus \Phi^\prime(0)\big)^2$, we get: 
        \begin{eqnarray*}
         \Delta & =  &  \big(1+  2  \proportion\sus  \Phi'(0)\big)^2 - 8 \proportion\sus \Phi'(0) +
         \big[1+  2  (\proportion\sus - \feedback)  \Phi'(0)\big]^2 -
           \big(1+  2  \proportion\sus  \Phi'(0)\big)^2                 \\
             & <  & (1 -   2  \proportion\sus  \Phi'(0))^2.   
         \end{eqnarray*}
Thus, $\sqrt{\Delta} < 1 - 2  \proportion\sus   \Phi' (0) $  and then:
         \begin{eqnarray*}
         \mu_1 & <   &  \frac{1 + 2  (\proportion\sus - \feedback)  \Phi' (0)   }{2} 
           + \frac{(1 - 2  \proportion\sus   \Phi' (0))}{2} \\
            & = & 1 - \feedback \Phi' (0) \\
             & < &  1. 
          \end{eqnarray*}   
Since the coefficient of $\mu$ in the characteristic polynomial is negative
we conclude that $ |\mu_2| < |\mu_1| < 1$ and hence the origin is locally asymptotically 
         stable.

\medskip
\noindent{\bf Case 3:  $\Delta > 0$,
        $  1 + 2  \proportion\sus \Phi' (0) < 
  2 \feedback \Phi' (0) < 2 + 4 \proportion\sus \Phi' (0) $
        and $2 \proportion\sus \Phi' (0) < 1$} (region 3 in Fig.~\ref{regions}).
       
By  simple algebraic manipulations, it follows that 
         $ -( 1 + 2  \proportion\sus \Phi' (0) ) <  ( 1 + 2  (\proportion\sus - \feedback)  \Phi' (0) )
                                         <  0. $      Therefore, as above, we obtain: 
          $\sqrt{\Delta} < 1 - 2  \proportion\sus   \Phi' (0) $  and consequently
         \begin{eqnarray*}
         \mu_2 & >  &  \frac{1 + 2  (\proportion\sus - \feedback)  \Phi' (0)   }{2} 
          - \frac{(1 - 2  \proportion\sus   \Phi' (0))}{2} \\
            & = & \frac{4 \proportion\sus \Phi' (0) - 2 \feedback  \Phi' (0)  }{2} \\
             & > & -1. 
          \end{eqnarray*}                                  
Since the coefficient of $\mu$ in the characteristic polynomial now is positive 
         we conclude that $ |\mu_1| < |\mu_2| < 1$ and hence the origin is locally asymptotically 
         stable.

\medskip   
\noindent{\bf Case 4:   $\Delta > 0$, $2 \feedback \Phi' (0) < 
        2 \proportion\sus \Phi' (0) -1  $
        and $2 \proportion\sus \Phi' (0) > 1$} (region 4 in Fig.~\ref{regions}).
       
Now, we have  $   2 (\proportion \sus -\feedback)  \Phi' (0) > 1  $ and, 
         therefore:
      \begin{eqnarray*}
       \mu_1 & >  & \frac{1  +  2 ( \proportion\sus - \feedback)  \Phi' (0)}{2}
 +  \sqrt{\Delta} \\
           & >  & \frac{1}{2}   + \frac{  2 ( \proportion\sus - \feedback) 
  \Phi' (0)}{2}  \\
           &  > & 1,   
       \end{eqnarray*}
proving that the origin is unstable in this case.

\medskip
  \noindent{\bf Case 5:  $\Delta > 0$,
                   $1+ 2 \proportion\sus  \Phi' (0)   <   2 \feedback \Phi' (0)  \leq 
        2+  4  \proportion\sus \Phi' (0)   $
        and $2 \proportion\sus \Phi' (0) > 1$} (region 5 in Fig.~\ref{regions}).

 Since $ \Delta > 0$, we must have either
 $  2 \feedback \Phi' (0) > \big(1 + \sqrt{ 2 \proportion\sus \Phi' (0)  }\, \big)^2 $ or
  $  2 \feedback \Phi' (0) < \big(1 - \sqrt{ 2 \proportion\sus \Phi' (0)  }\, \big)^2 $.
 But, since     $\big(1 - \sqrt{ 2 \proportion\sus \Phi' (0)  }\,\big)^2 < 1+ 2 \proportion\sus  \Phi' (0)
    <   2 \feedback \Phi' (0)$ then the second alternative is false.
 Thus:
    \begin{eqnarray*}
     2 \feedback \Phi' (0) &  >  &
 \big(1 + \sqrt{ 2 \proportion\sus \Phi' (0)  }\,\big)^2 \\
  &   = &
       1 + 2   \sqrt{ 2 \proportion\sus \Phi' (0)} + 
     2 \proportion\sus \Phi' (0)  \\
 &  > &  3 +  2 \proportion\sus \Phi' (0).  
  \end{eqnarray*}
We then obtain  $  2 (\proportion\sus - \feedback)  \Phi' (0) < -3$, from which
 it follows that
   \begin{eqnarray*}
       \mu_2 & <  & \frac{1  +  2 ( \proportion\sus - \feedback)  \Phi' (0)}{2}  \\
           & <  & -1,  
       \end{eqnarray*}
  proving instability.

\medskip
  \noindent{\bf Case 6:  $\Delta >0$,
                   $  2 \feedback \Phi' (0)  >    2+  4  \proportion\sus \Phi' (0)$} (region 6 in Fig.~\ref{regions}).
  
 In this case, since 
 $ 1  +  2   (\proportion J- \feedback  )\Phi' (0) < 
  1  +  2   \proportion\sus \Phi' (0) -2 - 4 \proportion\sus \Phi' (0) =  
   - 1  -  2   \proportion\sus \Phi' (0) < 0,$
  we have  
  $\Delta =  (1  +  2   (\proportion J- \feedback  )\Phi' (0))^2 - 8 \proportion\sus  > 
    (1  +  2    \proportion\sus \Phi' (0))^2 - 8 \proportion\sus  = 
    (1-   2   \proportion\sus   \Phi' (0))^2 $.
  
 If   $ 2 \proportion\sus  \Phi' (0)  < 1$,  then  $\sqrt{\Delta} > 1-2\proportion\sus \Phi' (0)$,  and 
  we obtain    
\begin{eqnarray*}
       \mu_2 & <  & \frac{1  +  2 ( \proportion\sus - \feedback)  \Phi' (0)}{2} 
     - \frac{1- 2   \proportion\sus  \Phi' (0)}{2}   \\
           & =  &   2   \proportion\sus  \Phi' (0) - \feedback  \Phi' (0)   \\
           & < &   2   \proportion\sus  \Phi' (0) - 1 -  2   \proportion\sus  \Phi' (0) \\
           & = & -1.
    \end{eqnarray*}
 
 If $ 2  \proportion\sus   \Phi' (0)> 1 $  then  $\sqrt{\Delta} > 2\proportion\sus \Phi' (0) - 1 $,
      and we obtain
   \begin{eqnarray*}
       \mu_2 & <  & \frac{1  +  2 ( \proportion\sus - \feedback)  \Phi' (0)}{2} 
     - \frac{- 1 +  2   \proportion\sus  \Phi' (0)}{2}   \\
           & =  &   1  - \feedback  \Phi' (0)   \\
           & < &   1 -  1 -  2   \proportion\sus  \Phi' (0) \\
           & < & -1.
    \end{eqnarray*}
Thus, $\mu_2<1$ in both case, and the instability follows.\hfill$\Box$

\subsection{Hopf bifurcation and periodic orbits}\label{ol} 
The occurrence of Hopf  bifurcation in the dynamical system (\ref{DYMA}) is the issue 
of  Theorem~\ref{limo} and Conjectures~\ref{noche} and \ref{checked} of the present 
section. All they are based on  
Theorem~\ref{hopf} below. It contains a description of the Hopf bifurcation phenomenon and 
provides conditions that are sufficient for it to occur. This bifurcation implies the emergence 
of periodic orbits in the dynamical system (\ref{DYMA}) for particular set of its 
parameters' values. 

\bigskip\noindent{\bf Theorem~\mishlabel{hopf} } (\cite{hale}, page 474, Poincar\'e-Hopf-Andronov 
theorem for maps).\newline
{\em    Let 
   \[ F: \R \times (\R)^2  \to \R^2; \quad (\eta,\mathbf{x}) \to F(\eta,\mathbf{x}) \]
   be a $\mathcal{C}^4 $ map depending on a real parameter $\eta$ satisfying the following conditions:
\begin{description}
\item{(i)} \ 
       $ F(\eta,0 ) = 0 $ for $\eta$ near some fixed $\eta_0$;
\item{(ii)} \   
       $DF(\eta,0 )$ (that is, $D_{\mathbf x}F(\eta,\mathbf{x})\big|_{\mathbf x=0}$) 
       has two non-real eigenvalues $\mu(\eta)$ and
      $  \bar{\mu}(\eta)$  for $\eta$ near  $\eta_0$, with $\big|\mu(\eta_0)\big|=1$;
\item{(iii)} \  
       $ \frac{d}{d \eta} \big|\mu(\eta)\big| >0 $ at $ \eta=\eta_0$;
\item{(iv)} \   
       $\mu^k (\eta_0) \neq 1$ for     $k=1,2,3,4.$  
\end{description} 
      Then there is a smooth $\eta$-dependent change of coordinates bringing $F$ into
       the form
        \[ F(\eta,\mathbf{x}) = \mathcal{F}( \eta,\mathbf{x}) + O(|| \mathbf{x}  ||^5  )  \] and there are
         smooth functions $A(\eta), B(\eta) \textrm{ and } \Omega(\eta) $ so that
         in polar coordinates the function  $\mathcal{F}( \eta,\mathbf{x})$  is given by
\begin{equation}\label{oprA}
              \left( \begin{array}{cc}
               r \\
                \theta 
                 \end{array} \right) \to 
                   \left( \begin{array}{cc}
               | \mu(\eta)| r - A(\eta) r^3 \\
                \theta + \Omega(\eta  ) + B(\eta) r^2
                 \end{array} \right) .
\end{equation}
           
If  $ A(\eta_0) >0, $   then there is a neighborhood $U$ of the  
origin  and a $\delta>0$ such that, for $|\eta- \eta_0| < \delta $ 
and $\mathbf{x} \in U$, the $\omega-$limit set of $\mathbf{x}$ is the origin 
if $\eta < \eta_0$, and belongs to a closed invariant  
$\mathcal{C}^1$ curve $\Gamma(\eta)$ encircling the origin if 
$\eta > \eta_0$. Furthermore  $\Gamma(\eta_0)=0$.\footnote{The phenomenon described 
                 in the last two sentences is called {\em Hopf bifurcation}. We add to this 
                 name the term {\em supercritical} in order to distinguish it from 
                 the case described in the next two sentences that will be called {\em the 
                 subcritical case}.} 

             If  $ A(\eta_0) < 0,$  then there is a neighborhood $U$ of  the 
             origin  and a $\delta>0$ such that, for $|\eta- \eta_0| < \delta $
              and $\mathbf{x}\in U$, the $\omega-$limit set of $\mathbf{x}$ is the origin 
               if $\eta > \eta_0$, and belongs to a closed invariant 
            $\mathcal{C}^1$ curve $\Gamma(\eta)$ encircling the origin if 
             $\eta < \eta_0$. Furthermore  $\Gamma(\eta_0)=0$.
}
             
\medskip\noindent{\bf Remark~\remlabel{hopfrem}.} \  We present here the method that 
we shall employ to calculate $A(\eta_0)$. The presentation follows \cite{hale}.
If the linear part of map $F$ at $\eta_0$
               is written in the Jordan canonical form,
\begin{equation}\label{nado} \mathbf{x} \to F(\eta_0, \mathbf{x}) = \left( \begin{array}{ll}
                 a  &  -b  \\
                  b  & a
                \end{array}  
                   \right)
                    \left( \begin{array}{l}
                  x_1 \\
                   x_2
                \end{array}  
                   \right) +
                    \left( \begin{array}{l}
                 G_1(x_1,x_2)   \\
                  G_2(x_1,x_2) 
                   \end{array}  
                   \right),\kern2em \mathbf{x}=(x_1, x_2),
\end{equation}
 where $a$ and $b$ relate to the eigenvalue $\mu(\eta)$ via the equation $\mu(\eta_0)=a+ib$, then for $A(\eta)$ defined in (\ref{oprA}), it holds that
\begin{equation}\label{mda} 
               A(\eta_0)= \textrm{Re} \left[\frac{(1-2 \mu(\eta_0)) \bar{\mu}(\eta_0)^2}{1 - \mu(\eta_0)} \xi_{11}\xi_{20 }  \right] 
                + \frac{1}{2} \left|\xi_{11}^2 \right| + \left|\xi_{02}^2 \right| - 
                 \textrm{Re}\left[\bar{\mu}(\eta_0) \xi_{21}\right],
\end{equation}
where (all the derivatives below are evaluated at $(x_1, x_2)=(0,0)$)
  \begin{description}
 \item $\xi_{20} = \frac{1}{8} \left\{
  (G_1)_{x_1x_1} -   (G_1)_{x_2x_2} + 2   (G_2)_{x_1x_2}  
 + i [  (G_2)_{x_1x_1} -   (G_2)_{x_2x_2}  -2  (G_1)_{x_1x_2}]  \right\}, $ 
 \item $\xi_{11} = \frac{1}{4} \left\{
  (G_1)_{x_1x_1} +  (G_1)_{x_2x_2} 
 + i [  (G_2)_{x_1x_1} +  (G_2)_{x_2x_2}  ]  \right\}, $ 
 \item $\xi_{02} = \frac{1}{8} \left\{
  (G_1)_{x_1x_1} -   (G_1)_{x_2x_2} - 2   (G_2)_{x_1x_2}  
 + i [  (G_2)_{x_1x_1} -   (G_2)_{x_2x_2}  -2  (G_1)_{x_1x_2}]  \right\}, $ 
\end{description}  
 and 
 \begin{eqnarray}
    \xi_{21}&  = & \frac{1}{16} \left\{
  (G_1)_{x_1x_1x_1} +   (G_1)_{x_1x_2x_2} +    (G_2)_{x_1x_1x_2} +  (G_2)_{x_2x_2x_2}
  \right.
\nonumber \\ 
  &  +  & \left. i [  (G_2)_{x_1x_1x_1} +   (G_2)_{x_1x_2x_2}  -  (G_1)_{x_1x_1x_2}
   -  (G_1)_{x_2x_2x_2}    ]  \right\}.\label{uau}   
 \end{eqnarray}

\bigskip
Theorem~\ref{limo} below reveals the occurrence of Hopf bifurcation in the dynamical system 
(\ref{DYMA}) when it satisfies certain additional assumptions. How general these assumptions are will
be discussed after the theorem's proof. The same discussion will expose at an intuitive lever 
the ideas behind the proof.

\bigskip\noindent{\bf Theorem~\mishlabel{limo}}.\newline
{\em  Suppose  the dynamical system (\ref{DYMA})  satisfies the following assumptions:
\begin{description}
\item{{\bf Assumption \assone}.} The  parameters $\feedback, \proportion, \sus $
  and the distribution 
function $\Phi(\cdot)$  are
 all functions of a single real variable $\eta$ that are defined on some 
 (nonempty) interval  $[\eta_\ini, \, \eta_\fini]$; the functions $\feedback_\eta, \proportion_\eta, \sus_\eta$ and the distribution functions  $\{\Phi_\eta(\cdot)\}$ 
 are  smooth enough to ensure that 
  $\Psi$, the map that generates the system (\ref{DYMA}), is $\mathcal{C}^4$ 
 in $(\price_{n-1},\demand_{n-1} )$  and $\eta$. 
\item{{\bf Assumption \asstwo}.} The functions  $\proportion_\eta, \sus_\eta$ and 
$\Phi^\prime_\eta(0)$ (the latter means the value of  $\partial\Phi_\eta(x)/\partial x$ at $x=0$) 
are all increasing in $\eta$ on the domain $[\eta_\ini, \, \eta_\fini]$.
\item{{\bf Assumption \assthree}.} 
$(1 + 2 (\proportion_\eta\sus_\eta - \feedback_\eta) \Phi_\eta'(0))^2 - 8 \proportion_\eta\sus_\eta \Phi'_\eta(0))<0,  \hbox{ for each }\eta\in [\eta_\ini, \, \eta_\fini]$.
\item{{\bf Assumption \assfour}.}  There exists a unique $\eta_0\in [\eta_\ini, \, \eta_\fini]$ such that 
$2\proportion_{\eta_0}\sus_{\eta_0}\Phi^\prime_{\eta_0}(0)=1$, and $\eta_0\not=\eta_\ini$, $\eta_0\not= \eta_\fini$.
\end{description}
Suppose in addition that
\begin{equation}\label{apred}
(a)\,\,\feedback_{\eta_0} = \frac{1}{2},\kern1em (b)\,\, \Phi'_{\eta_0}(0)=1,\kern1em (c)\,\, \Phi^{\prime\prime}_{\eta_0}(0)=0,\kern1em (d)\Phi^{\prime\prime\prime}_{\eta_0}(0)< 0.
\end{equation}

Then the system undergoes the supercritical Hopf bifurcation when $\eta$ passes 
through $\eta_0$.
}

\medskip\noindent{\em Proof.} \ In virtue of Thm.~\ref{hopf} and Remark~\ref{hopfrem}, 
in order to prove the present theorem, it is sufficient to show that the mapping $\Psi$ that 
generates the dynamical system (\ref{DYMA}) satisfies the conditions (i) -- (iv) 
of Thm.~\ref{hopf} and that $A(\eta_0)$ is a positive number. We do so in Steps~1 -- 5 below.

{\em Step 1:} \ Since Theorem~\ref{propdyns} ensures that $(0,0)$ is an equilibrium point 
of (\ref{DYMA}) for any values of the parameters of this dynamical system, then 
$\Psi$ satisfies (i) for each $\eta\in [\eta_\ini, \, \eta_\fini]$. 
 
{\em Step 2:} \ We start recalling facts and results from the proof of Thm.~\ref{propdyns} 
that we shall need below: (a) the Jacobian matrix of the mapping $\Psi$ at $(0,0)$ (i.e., 
$DF(\eta,0 )$, in the notations of Thm.~\ref{hopf}) was 
calculated and the result is presented in (\ref{jacob}); (b) the matrix' eigenvalues were 
calculated and their expressions are presented in (\ref{rutik}); (c) it was proved that 
these eigenvalues are non-real numbers provided $\Delta<0$, with $\Delta$ being defined by
(\ref{Delta}); (d) the modulus of each eigenvalue was found to be equal to 
$\sqrt{2\proportion \sus\Phi^\prime(0)}$, provided $\Delta\leq 0$. 

Now, since Ass.~\assthree\ ensures that $\Delta<0$ for every $\eta$, then (c) above 
implies that $\Psi$ satisfies the first part of the condition (ii). As for the second part 
of this condition, it is implied by Assumptions~\assthree and \assfour\ and the the fact 
(d). 

{\em Step 3:} \ The fact (d) from the list of Step~2 and  Assumptions~\asstwo\ and \assthree\ 
imply straightforwardly that  the condition (iii) is satisfied by the eigenvalues of $\Psi$ 
for every $\eta$, and in particular, for $\eta_0$.

{\em Step 4:} \ {}From the fact (b) of the list of Step~2 and from Ass.~\assthree, it follows that 
\begin{eqnarray}
\mu(\eta_0) &= &\frac{1 + 2( \proportion_{\eta_0} \sus_{\eta_0}-\feedback_{\eta_0}) \Phi'_{\eta_0}(0)}{2}\label{kore}\\
 &&  + i \frac{\sqrt{8 \proportion_{\eta_0} \sus_{\eta_0} \Phi'_{\eta_0}(0)) -
  (1 + 2( \proportion_{\eta_0} \sus_{\eta_0}-\feedback_{\eta_0}) \Phi'_{\eta_0}(0)
      )^2 }}{2 }.\nonumber 
\end{eqnarray}
Applying Ass.~\assfour\ and constraints (\ref{apred}) to (\ref{kore}), we get that   
$\mu(\eta_0)=\frac{1}{2} + i \frac{\sqrt{3}}{2}$. {}From this, $\mu^k(\eta_0)\not=1, k=1,2,3,4$, 
proving hence the validity of condition (iv).

{\em Step 5:} \ In this step, we shall prove that $A(\eta_0)>0$. In the calculations that follow, 
all functions that depend on $\eta$ will be evaluated at $\eta_0$. This allows us to omit 
$\eta_0$ in the notations throughout the proof, namely, we shall write $A$, $\feedback, \proportion, 
\sus$ and $\Phi^\prime(0)$ for, respectively, $A(\eta_0)$, $\feedback_{\eta_0}, \proportion_{\eta_0}, 
\sus_{\eta_0}$ and $\Phi_{\eta_0}^\prime(0)$. 

Since we intend to use the expression (\ref{mda}) then we need to get a (\ref{nado})-like 
form of the map $\Psi$. We shall get it from the Jacobian matrix  at $(0,0)$ of 
the map $\Psi$ that we've calculated in the proof of Thm.~\ref{propdyns} (the matrix is 
presented in (\ref{jacob})). Since  Ass.~\assfour\ 
imposes that $2\proportion\sus\Phi^\prime(0)=1$ then this matrix acquires (recall, at $\eta_0$) 
the following form: 
\[
M=\left( \begin{array}{cc}
 1 & \feedback \\
 -2 \Phi^\prime(0) &   1-2 \feedback\Phi^\prime(0) 
  \end{array} \right).
\]
It is then easy to check that the matrix 
 \[ 
 { P} = \left( \begin{array}{cc}
 1 & 0 \\
 -\Phi^\prime(0) & \frac{\sqrt{1- \varepsilon^2 }}{\feedback}
   \end{array} \right), \quad \hbox{with}\quad \varepsilon:= 1- \feedback \Phi^\prime(0),
\]
puts $M$ in the Jordan canonical form: 
$P^{-1}MP=\left(  \begin{array}{cc}
    \varepsilon &  \sqrt{1-\varepsilon^2  } \\
     -\sqrt{1-\varepsilon^2  } & \varepsilon
\end{array}  \right)$. 
Therefore,  in the new variables $ \left( \begin{array}{c}
    \tilde\price_n    \\ \tilde \demand_n   \end{array} \right) :=
   { P}^{-1}    \left( \begin{array}{c}
       \price_n \\ \demand_n   \end{array} \right)$, $n\in \mathbb{N}$,
the map   $\Psi$ acquires a  (\ref{nado})-like form:
\begin{eqnarray*}
\left(
  \begin{array}{cc} \tprice_n \\ \tdemand_n \end{array} \right) &  = &
 P^{-1} M P  
 \left(\begin{array}{cc} \tprice_{n-1} \\  \tdemand_{n-1} \end{array} \right)
\,\,   + P^{-1} \left(\Psi-M\right) P
 \left(\begin{array}{cc} \tprice_{n-1} \\  \tdemand_{n-1} \end{array} \right)\\
&=& \left(
   \begin{array}{cc}
    \varepsilon &  \sqrt{1-\varepsilon^2  } \\
     -\sqrt{1-\varepsilon^2  } & \varepsilon
\end{array} \right)   
  \left(
  \begin{array}{cc}
     \tprice_{n-1} \\
     \tdemand_{n-1}
   \end{array}     \right)
   +  \left(
   \begin{array}{c}
      G_1 (\tprice_{n-1},\, \tdemand_{n-1})\\
       G_2( \tprice_{n-1}, \tdemand_{n-1})  
   \end{array}     
     \right),
\end{eqnarray*}
where
\begin{eqnarray*} G_1(x_1, x_2)&=&0,\\
G_2(x_1, x_2) &=& 
   \frac{\feedback}{\sqrt{1-\varepsilon^2}}\Bigl\{ 
 \proportion \left[1-2 \Phi\bigl(K_1 x_1 + K_2 x_2\bigr)\right] +
  (1-\proportion)  \left[1-2 \Phi\bigl(K_3 x_1 + K_4 x_2\bigr)\right] \\
& &\kern1em + \Phi^\prime(0)\left(3-2\proportion\Phi^\prime(0)\right)x_1-
\feedback^{-1}\left(1-2\proportion \Phi^\prime(0)\right)\sqrt{1-\varepsilon^2}
x_2\Bigr\},
\end{eqnarray*}
and where
\[\begin{array}{rclrcl}
  K_1 & = & 1 -(\feedback-\sus)\Phi^\prime(0), \kern2em  
  & K_2 & = & {\feedback}^{-1}\sqrt{1-\varepsilon^2}(\feedback- \sus), \\
  K_3 &= & 1-\feedback \Phi^\prime(0), &  K_4 & = & \sqrt{1-\varepsilon^2}.
\end{array}\]

{}From the expressions for $G_1$ and $G_2$, we shall now obtain  $\xi_{\ell r}$'s 
following the formulas from Remark~\ref{hopfrem}.  Differentiating $G_2$  we get:
\[
\frac{\partial^2}{\partial x_1^2}G_2(x_1, x_2) = -\frac{2\feedback}{\sqrt{1-\varepsilon^2}}
\left(\proportion K_1^2\Phi^{\prime\prime}(K_1x_1+K_2x_2)+
                          (1-\proportion) K_3^2 \Phi^{\prime\prime}(K_3x_1+K_4x_2) \right).
\]
{}From this result, in virtue  of the assumption (\ref{apred}-c), we conclude that 
$\frac{\partial^2}{\partial x_1^2}G_2(0,0)=0$. By similar arguments, we conclude that 
all second order derivatives of $G_2$ vanish at $(x_1, x_2)=(0, 0)$. These facts and 
the fact that $G_1\equiv 0$ imply that $\xi_{20} = \xi_{11} = \xi_{02} = 0$. In order
to find $\xi_{21}$, we calculate third order derivatives of $G_2$:
\begin{eqnarray*}   
\frac{\partial^3}{\partial x_1^3} G_2(0,0) & = &
 - \frac{2\feedback\Phi'''(0)}{\sqrt{1-\varepsilon^2}} \left(\proportion K_1^3 + 
   (1-\proportion) K_3^3\right), \\
\frac{\partial^3}{\partial x_1^2\partial x_2} G_2(0,0) & = &
 - \frac{2\feedback\Phi'''(0)}{\sqrt{1-\varepsilon^2}} \left( \proportion K_1^2 K_2 + 
   (1-\proportion) K_3^2 K_4\right), \\
\frac{\partial^3}{\partial x_1\partial x_2^2} G_2(0,0) &=&  
 - \frac{2\feedback\Phi'''(0)}{\sqrt{1-\varepsilon^2}} \left( 
 \proportion K_1 K_2^2 +    (1-\proportion) K_3 K_4^2\right), \\
\frac{\partial^3}{\partial x_2^3} G_2(0,0) &=&  
 - \frac{2\feedback\Phi'''(0)}{\sqrt{1-\varepsilon^2}} \left( 
 \proportion  K_2^3 +    (1-\proportion)  K_4^3\right). 
 \end{eqnarray*}
Now, from the relations  $2\proportion \sus\Phi^\prime(0)=1$ (that holds because of 
Ass.~\assfour) and $\Phi^\prime(0)=1$ (ensured by (\ref{apred}-b)) we get that 
$\proportion=(2\sus)^{-1}$. This allows us to substitute $\proportion$ by $(2\sus)^{-1}$ 
in the expressions for the derivatives of $G_2$. We also substitute there $\feedback$ 
and $\Phi^\prime(0)$ by $\frac{1}{2}$ and $1$, respectively (these 
substitutions are justified by (\ref{apred}-a, b)). 
The resulting simplified expressions for the 
derivatives lead, via the formula (\ref{uau}), to the following: 
\[\begin{array}{rcl}
\xi_{21} &=& \frac{\Phi'''(0)}{64 \sqrt{3}}[ (-20 \sus^2 +
  (13 \sqrt{3}-4) \sus + 3 \sqrt{3}-10 )\\
   &&+ i   (-10 \sus^2 +
  (4 \sqrt{3}-7) \sus + 2 \sqrt{3}-5 )  ].\end{array}
\]
Plugging in (\ref{mda}) the expressions for $\xi_{\ell r}$'s obtained above, we 
finally get that
\[    A  =  - \textrm{Re}(\bar{\mu} \xi_{21})
    =  - \Phi'''(0) \frac{2+\sqrt{3}}{64 \sqrt{3}}
    [5 \sus^2 + (1- 2 \sqrt{3})\sus +1].
\]
This expression for $A$ and the assumption (\ref{apred}-d) ensure that $A>0$ for every 
$\sus$. This completes Step 5 and the proof of Theorem~\ref{limo}.\hfill$\Box$

\bigskip We proceed with the discussion of how general the Hopf bifurcation phenomenon is 
for the dynamical system (\ref{DYMA}). 

In order to discuss Hopf bifurcation in a dynamical system with the help of 
Theorem~\ref{hopf}, the minimal necessary condition is that all system's parameters 
be expressed as functions of a unique variable. In the framework of Theorem~\ref{limo}, 
this condition is ensured by Assumption~\assone. Let us accept it now and let us 
consider then two functions: $u_\eta:=2\proportion_\eta \sus_\eta\Phi^\prime_\eta(0)$ 
and $w_\eta:= 2\feedback_\eta \Phi^\prime_\eta(0)$. We recall from the proof of 
Thm.~\ref{propdyns} that {\sf (a)} the quantity denoted there by $\Delta$ can be expresses as a 
function of solely $u$ and $w$, {\sf (b)} when $\Delta(u, w)<0$ then the eigenvalues of 
the linearization of $\Psi$ at its fixed point $(0,0)$ are non-real numbers, and {\sf (c)} when 
$u=1$ then the modulus of each eigenvalue is $1$. The facts {\sf (a,b,c,)} ensure that 
if the parameter $\eta$ introduced in Assumption~\assone\ is such that $\Delta(u_\eta, w_\eta)<0$, 
for all $\eta$, and if $u_{\eta_0}=1$, for some $\eta_0$, then the property (ii) of 
Thm.~\ref{hopf} is satisfied. 
These ``if''\ conditions are provided by Assumptions~\assthree\ and \assfour. 
As for the Ass.~\asstwo\, ensures the validity of (iii) of Thm.~\ref{hopf}.  
For this, Assumption~\asstwo\ may be not the minimal sufficient condition, but 
we did not search for such. 
Turning our attention to the condition (i) of Thm.~\ref{hopf}, we easily see that it
is valid for $\Psi$ due to Thm.\ref{propdyns}(a) without the necessity for any additional assumptions.

We thus have showed that Assumptions~\assone -- \assfour\  ensure the 
validity of conditions (i) -- (iii) of Thm.~\ref{hopf}. It seems to us that 
these assumptions are also sufficient for the validity of the 
condition (iv) of Thm.~\ref{hopf} and the inequality $A(\eta_0)\not = 0$.
Accordingly, we formulate the following:

\bigskip\noindent{\bf Conjecture~\conjlabel{noche}.}\newline
{\em The dynamical system (\ref{DYMA}) exhibits periodic orbits, when the initial point
is sufficiently close to $(0,0)$ and when the parameters are such that the corresponding 
point $(u,w)$ is sufficiently close to the right of the interval 
$\{(u,w)\,:\, u=1,\, 1<w<4\}$;\footnote{If the point is close to this interval on its left, 
then $(0,0)$ is locally asymptotically stable equilibrium of the dynamical system (\ref{DYMA}) 
and hence the periodic orbits described in this statement cannot occur.}
precisely to state, there is a neighborhood $U$  of $(0,0)$
such that the $\omega$-limit of any orbit of (\ref{DYMA}) with initial condition
in $U$ belongs to a closed $\mathcal{C}^1$ curve encircling $(0,0)$, provided  
$2\proportion \sus\Phi^\prime(0)$ is slightly bigger  than $1$ and 
$\feedback \Phi^\prime(0)$ is between $0$ and $2$.}

\bigskip We did not  pursue proofs in the degree of generality that would allow us to justify 
rigorously the generic property formulated in this conjecture. Rather, we considered two
particular cases. The first case is presented in Thm.~\ref{limo}. There, we assumed 
an additional constraint (\ref{apred}) that helped a lot to simplify the calculations 
needed to establish the condition (iv) and the inequality $A(\eta_0)>0$. The 
second case is presented in the statement below. This case attracted our attention 
because it arises in applications of our mathematical 
study of the dynamical system (\ref{DYMA}). The reason for this is explained in 
Section~\ref{application}.

\bigskip\noindent{\bf Conjecture~\conjlabel{checked}.}\newline
{\em  If the parameters of the dynamical system (\ref{DYMA}) satisfy 
Assumptions~\assone\ -- \assfour\ and $\Phi$ is Normal Distribution 
with null mean then the system undergoes the supercritical Hopf bifurcation 
when $\eta$ passes increasingly through $\eta_0$.} 

\bigskip In the case of Conjecture~\ref{checked}, the analytic verification
of the condition (iv) and the inequality  $A(\eta_0)\not=0$ turned to be an extremely 
tedious task. To carry out this task, we resorted to numeric methods: we verified the condition 
and the inequality numerically for a grid of parameter  values.  

\subsection{Global asymptotic stability for small values of the parameter $\lambda$}\label{gogo}
 
\bigskip\noindent{\bf Theorem~\mishlabel{dyns}} \ (a sufficient condition for global asymptotic stability).\newline 
{\em  Suppose $ 2 \proportion \sus \Phi'(0) < 1$. Then, there exists a 
positive number $\feedback_c$ such that if  $\feedback \in (0, \feedback_c]$  
then the equilibrium $(0,0)$ of the dynamical system (\ref{DYMA}) is  globally  
asymptotically stable. (The numeric value of $\feedback_c$ will be 
specified in the proof by (\ref{cico}).)}

\medskip\noindent{\bf Remark~\remlabel{keshe}.} \ On comparing  the assertions of 
Thms.~\ref{propdyns} and~\ref{dyns} one concludes that the assumptions of 
Thm.~\ref{dyns} must be a sub-case of the assumptions of Thm.~\ref{propdyns}. 
This inclusion can be established by simple calculations grounded on 
(\ref{cico}-(\feddois)), one of the conditions that will determine the value 
of $\feedback_c$. We shall omit these calculations.

\bigskip We now proceed with the argument that will lead to the proof of
Thm.~\ref{dyns}. The arguments employ essentially the assumption (\ref{parcon})
that ensures that $\sus>0$ and $\proportion>0$ and thus, allows us to divide
by $\sus$ and $\proportion$. Actually, the theorem's assertion remains true
even when either $\sus$ and $\proportion$ is equal to $0$, but this case 
requires a specific argument which will not be presented here. 

Upon the substitution of $\price_{n-1}+\feedback \demand_{n-1}$ by $\price_n$ in 
the second equation of the dynamical system (\ref{DYMA}), it acquires the form 
that makes it clear that the passage from $(\price_{n-1}, \demand_{n-1})$ to
$(\price_{n}, \demand_{n})$ can be considered as a two-step procedure: first, 
$\price_n$ is obtained from $\price_{n-1}$ and  $\demand_{n-1}$ via the first 
equation of the system, and then, $\demand_n$ is obtained from $\demand_{n-1}$ 
and $\price_n$ via the second equation. The second equation, or better to say, 
its non-linearity is the principal cause for the difficulty of studying the evolution
of the dynamical system (\ref{DYMA}). In the present proof, we overcome the 
difficulty by getting convenient estimates of two quantities (to be defined in  
(\ref{kabe}) and (\ref{mabe}) below) that relate to this equation and characterize 
its features that are important for our analysis. The presentation of these estimates 
becomes more transparent -- in our opinion --, when we express the equation through an 
appropriate one-parameter $\mathbb{R}\rightarrow \mathbb{R}$ mapping in which $\price_n$ 
is the parameter and $\demand_{n-1}$ is the  argument. In order to do so, for the 
``parameter'' $\price\in \mathbb{R}$, we define $g_\price (\cdot)\,:\, (-\infty, \infty)
\rightarrow [-1,1]$ as follows:\footnote{Although the argument of $g_\price(\cdot)$ cannot 
     exceed $1$ in modulus -- since it has been defined as the mathematical analogue 
     of population excess demand --, it turns out to be convenient to extend the 
     function's argument domain to the whole $\mathbb{R}$.}
\begin{equation}\label{naga}
g_\price (\demand) = 
  \proportion [1-2\Phi \big( \price  - \sus \demand\big )] + 
  (1-\proportion) [1-2\Phi \big(\price \big)],\kern1em d\in (-\infty, \infty).
\end{equation}
Now, we can re-write (\ref{DYMA}) in the desired form:
\begin{eqnarray}
\price_n&=&\price_{n-1}+\feedback \demand_{n-1}\label{clan}\\
\demand_n&=&g_{\price_{n}}\left(\demand_{n-1}\right)\label{blan}
\end{eqnarray}

From now  on, we shall use the shorthand notation:
\begin{equation}\label{debe}
b:= 2 \proportion \sus \Phi'(0);
\end{equation}
note that $b>0$ (because of (\ref{phicon}-(\diftza), (\mondir))) and $b<1$ 
(because of the theorem's assumption).

We start with a list of basic properties of $g_\price(\cdot)$. For each $\price$, 
\begin{eqnarray}
\circ&&g_\price(\cdot) \hbox{ is everywhere differentiable function and } g^\prime_\price(\demand)=2\proportion\sus\Phi^\prime(\price -\sus\demand);\label{edive}\\
\circ&&\sup_{d\in\mathbb{R}}\{g^\prime_\price(d)\}=g_\price^\prime\left(\frac{\price}{\sus}\right)=2\proportion\sus\Phi^\prime(0)=b\in (0,1);\label{gadiv}\\
\circ&&\hbox{in particular, as a consequence of (\ref{edive}) and (\ref{gadiv}), } g_\price(\cdot)\hbox{ is a contraction: }\nonumber\\ &&\kern6em|g_\price(x)-g_\price(y)|\leq b|x-y|<|x-y|;\label{squez}\\
\circ&& g_\price(\cdot)\hbox{ a monotone strictly increasing function};
\label{moveh}\\
\circ&&\hbox{the graph of }g_\price(\cdot) \hbox{ passes through the point }\nonumber\\  
&&\kern4em\left(\demand^\ast, y^\ast\right)\hbox{ with } \demand^\ast=\frac{\price}{\sus}\hbox{ and }
y^\ast=(1-\proportion)[1-2\Phi(\price)]\label{tupo}\\
&&\kern1em\hbox{that satisfy: }\left\{\begin{array}{l}\price\geq\, (>)\, 0\\ \price\leq \,(<)\, 0\end{array}\right.\,\,
 \Rightarrow\,\, \left\{\begin{array}{l}\demand^\ast\geq \, (>)\, 0,\, y^\ast\leq\,(<)\, 0,\\
 \demand^\ast\leq \, (<)\, 0,\,
y^\ast\geq \, (>)\, 0.\end{array}\right. \nonumber%
\end{eqnarray}
All these properties stem from (\ref{naga}) in combination with the properties  of 
$\Phi$ assumed in (\ref{phicon}); exactly to state, we used (\ref{phicon}-\diftza) 
for (\ref{edive}), (\ref{phicon}-\maxdir) for (\ref{gadiv}), (\ref{phicon}-\mondir) 
for (\ref{moveh}), and (\ref{phicon}-\eqibda) for (\ref{tupo}).  

\begin{figure}[h] 
\begin{center}
\begin{picture}(400,350)
\put(0,250){\vector(1,0){400}} \put(390,230){$d$}
\put(200,40){\vector(0,1){310}}\put(210,330){$g_p(d)$} 
\put(0,50){\line(1,1){300}} 
\qbezier[25](180,250)(175,240)(185,235)\put(150,230){$45^{\circ}$}
\put(0,50){\line(3,2){400}} 
\put(240,210){\circle*{3}}  
\put(300,250){\circle*{3}}  
\put(290,270){$\mathcal{R}_{{\rm lower}}$}
\put(0,50){\circle*{3}}     
\multiput(1,50)(0,2){101}{\circle*{.1}}
\put(10,270){$\tilde{\mathcal{D}}(p)$} 
\put(390,260){\line(-3,-1){180}} 
\put(360,250){\circle*{3}}  
\put(370,270){$\mathcal{R}_{{\rm upper}}$}
\thicklines
\qbezier(240,210)(300,250)(400,270)
\qbezier(0,120)(150,150)(240,210)  
\thinlines
\put(93,143){\circle*{3}}  
\multiput(93,143)(0,2){54}{\circle*{.1}}\put(93,270){$\mathcal{D}(p)$}
\put(325,250){\circle*{3}}  
\put(330,270){$\mathcal{R}(p)$}
\qbezier[15](270,230)(290,220)(285,210)\put(300,220){$\vartheta$}\put(250,170){$\tan\vartheta=b$\kern3em$\tan\xi=a$}
\qbezier[30](330,240)(345,225)(340,210)\put(345,220){$\xi$}
\multiput(200,210)(2,0){100}{\circle*{.1}}
\multiput(240,210)(0,2){21}{\circle*{.1}}
\put(74,210){$(1-\alpha)\left[1-2\Phi(p)\right]=:y^\ast$}
\put(240,265){$d^\ast:=\frac{p}{J}$}
\put(160,130){$\alpha\left[1-2\Phi(p-Jd)\right]+(1-\alpha)\left[1-2\Phi(p)\right]$}
\put(155,135){\vector(-3,2){30}}
\put(110,90){the tangent line at the point $(d^\ast, y^\ast)$}
\put(105,95){\vector(-3,2){18}}
\end{picture}
\end{center}
\centering
\caption{Illustration to the definition of the quantities $\mathcal{D}(p)$ and $\mathcal{R}(p)$ and the arguments that employ these quantities and their properties.}\label{yA}
\end{figure} 

We now define the quantities $\mathcal{D}(p)$ and $\mathcal{R}(p)$ 
that will play central role in our arguments:
\begin{eqnarray}
\mathcal{D}(p)&:=&\hbox{the solution of }g_\price(\demand)=\demand,\hbox{ hence, }
g_\price(\mathcal{D}(p))=\mathcal{D}(p),\label{kabe}\\
\mathcal{R}(p)&:=&\left\{\begin{array}{l}
-\infty,\hbox{ when }g_\price(\demand)>0,\,\, \forall\demand,\\
+\infty,\hbox{ when }g_\price(\demand)<0,\,\, \forall\demand,\\
\hbox{the solution of }g_\price(\demand)=0,\hbox{ otherwise, and hence, }g_\price(\mathcal{R}(p))=0.
\end{array}\right.\label{mabe}
\end{eqnarray}
Figure~\ref{yA} helps  the  visualization of  these definitions. We observe that the properties 
(\ref{gadiv}) and (\ref{moveh}) ensure that $\mathcal{D}(p)$ and $\mathcal{R}(p)$ are well and 
uniquely defined.

We shall frequently use the following properties of $\mathcal{D}(p)$ and $\mathcal{R}(p)$: 
\begin{eqnarray}
\circ&&\hbox{ if }\price\geq 0\hbox{ then }\mathcal{R}(\price)\geq\frac{\price}{\sus}\geq 0,\hbox{ and if }\price\leq 0\hbox{ then }\mathcal{R}(\price)\leq \frac{\price}{\sus}\leq 0;\label{tvah}\\
\circ&&\hbox{ if }\price\leq \price^\prime \hbox{ then }\mathcal{D}(\price)\geq \mathcal{D}(\price^\prime);\label{juti}\\
\circ&&\hbox{ if }\price\geq 0\hbox{ then }\mathcal{D}(\price)\leq 0,\hbox{ and if }\price\leq 0\hbox{ then }\mathcal{D}(\price)\geq 0.\label{pozd} 
\end{eqnarray} 
(\ref{tvah}) follows directly  from (\ref{moveh}) and (\ref{tupo}). To prove (\ref{juti}), we 
argue as follows: First, from the monotonicity of $\Phi$ and the definition (\ref{naga}), we get:
\begin{equation}\label{monip}
\circ\kern2em g_{\price}(\cdot)\hbox{ is decreasing in }\price, \hbox{i.e., }g_{\price}(\demand)>g_{\price^\prime}(\demand)\, \forall\, \demand,\hbox{ when }\price<\price^\prime.
\end{equation}
Then, from this property and the fact that $|g_\price^\prime(\demand)|<1$, we conclude 
that the first coordinate  of the intersection point of  the $45^\circ$ line in the plane 
with $g_{\price}(\cdot)$ decreases with the increase of $\price$. This conclusion is exactly 
the property (\ref{juti}). As for (\ref{pozd}), it follows from (\ref{juti}) because 
$\mathcal{D}(0)=0$ (this equality holds because $g_0(\cdot)$ passes through 
$(0,0)$ as ensured by (\ref{phicon}-\eqibda) and (\ref{naga})).  

There are two more properties of $\mathcal{D}(\price)$ and $\mathcal{R}(\price)$ that we 
shall need for the proof of Theorem~\ref{dyns}. These are provided by Lemmas~\ref{estAh} 
and \ref{estrh} below.

\bigskip\noindent{\bf Lemma~\ammellabel{estAh}.}\ \ {\em
   If  $\price \geq 0$   then     
   $\mathcal{D}(\price) \geq - \frac{b}{\sus \proportion (1-b)} \price$, and if
   $\price \leq 0$ then  $ \mathcal{D}(\price)  \leq - \frac{b}{\sus \proportion (1-b)} \price$.
}    

\medskip\noindent{\bf Proof. } \ We prove the case $\price\geq 0$; in the case $\price\leq 0$, 
the proof is similar. 

We shall need the following estimate:
\begin{equation}\label{tena}
1-2\Phi(\price) \stackrel{\rm (\ref{phicon}\,-\,(\diftza),\, (\eqibda))}{=}
-2\int^{\price}_0\Phi^\prime (x)dx\stackrel{\rm (\ref{phicon}-(\maxdir))}{\geq}  -2\Phi^\prime(0)\price \stackrel{\rm(\ref{debe})}{=}-\frac{b\price}{\proportion\sus}.
\end{equation}

Now\footnote{The following constructions are illustrated in Figure~\ref{yA}.}, 
we construct the tangent line to $g_\price(\cdot)$ at the point 
$(\demand^\ast, y^\ast)$. From (\ref{edive}), (\ref{gadiv}) and (\ref{tupo}), we have 
that the slope of the line is $b$ and the graph of $g_\price(\cdot)$ to the left of 
$(\demand^\ast, y^\ast)$ lies above this line. We also have that $\demand^\ast\geq 0$ 
and $y^\ast\leq 0$ (because of (\ref{tupo}) and the assumption $\price\geq 0$). Hence, 
$\mathcal{D}(p)$, the coordinate of the intersection of $g_\price(\cdot)$ 
with the diagonal, is larger than $\tilde{\mathcal{D}}(\price)$, the coordinate of 
the intersection of the tangent line with the diagonal. By simple calculations, we get that
$\tilde{\mathcal{D}}(\price)=(1-b)^{-1}\left\{ (1-\proportion)\bigl[1-2\Phi(p)\bigr]-
(b\price)/\sus\right\}$. Combining this with the inequality  $\mathcal{D}(\price) \geq \tilde{\mathcal{D}}(\price)$ just obtained and with the estimate (\ref{tena}), we 
complete the proof: $\mathcal{D}(\price)\geq \tilde{\mathcal{D}}(\price)  \geq     
\frac{ - \frac{1-\proportion}{\proportion} b \price  -b\price}{\sus(1-b)} = 
- \frac{b}{\sus \proportion (1-b)}\price$.\hfill $\Box$

\bigskip\noindent{\bf Lemma~\ammellabel{estrh}.}\ \ {\em Let ($a$ defined 
below might have been any other number from  $(0,\,b)$):
\begin{equation}\label{aaa}
a:=b/2.
\end{equation}
Let next $x_\ell<0$, $x_r>0$ be such that:
\begin{equation}\label{vera}
2\proportion\sus \Phi^\prime(x)\geq a,\,\,\,\forall x\in [x_\ell, x_r]
\end{equation} 
(the existence of $x_\ell$ and $x_r$ follows from the properties (\ref{phicon}-\diftza, \maxdir, \mondir) 
of $\Phi$; in fact, $x_\ell$ and $x_r$ are the negative and the positive solutions of the equation $2\proportion\sus \Phi^\prime(x)=a$). Finally: 
\begin{equation}\label{hhh}
h:=\min\left\{\frac{\proportion}{1-\proportion}\frac{a}{b},\,\, 1\right\}
\,\times\,{\min\left\{|x_\ell|, x_r\right\}}.
\end{equation}

Then, $\mathcal{R}(\cdot)$ satisfies the following properties:}\newline 
{\bf (a)} {\em If $\price\geq 0$ then
     $\mathcal{R}(\price)\geq  \frac{1}{\proportion \sus} \frac{a}{b} \price $, 
     and if $\price\leq 0$ then $\mathcal{R}(\price) \leq   
       - \frac{1}{\proportion \sus} \frac{a}{b} \price$.} \newline
{\bf (b)} 
{\em If $0\leq \price\leq h$ then
     $\mathcal{R}(\price) \leq   
        \frac{1}{\proportion \sus} \frac{b}{a} \price $, and if $-h\leq \price\leq 0$ then 
$\mathcal{R}(\price)\geq  -\frac{1}{\proportion\sus} \frac{b}{a} \price $. 
}

\medskip\noindent {\bf Proof.} \ We shall prove both (a) and (b) for $\price \geq 0$. For 
$\price \leq 0$, the proofs are analogous. 

Due to (\ref{hhh}), if  $0\leq \price\leq h$ then $x\in [-\price, 0]$ implies $x\in [x_\ell, 0]$.
This implication and (\ref{vera}) ensure that $\Phi^\prime(x)\geq (2\proportion \sus)^{-1}a$, in case
$x\in [-\price, 0]$ and  $\price\in [0,h]$. We use this fact to derive the following inequality:  
\begin{equation}\label{estF}
1-2\Phi(\price)\stackrel{\rm (\ref{phicon}\,-\,(\diftza),\, (\eqibda))}{=}
-2\int^{\price}_0\Phi^\prime (x)dx
{\leq} -\frac{a}{\proportion\sus}\price, \kern1em\hbox{ for all }\price\in [0,h].
\end{equation}

Recall now from (\ref{tupo}) that the graph of $g_\price$ passes trough the point 
$(\demand^\ast, y^\ast)$. We construct two straight lines passing through this point 
and having  slopes $b$ and $a$, and  we denote by $\mathcal{R}_{\rm lower}$ and 
$\mathcal{R}_{\rm upper}$  the coordinates of their respective intersections with 
the $\demand$ axis (see Figure~\ref{yA}). Directly from this construction, we
have that $\mathcal{R}_{\rm lower} =  \demand^\ast-\frac{y^\ast}{b}$ and  
$\mathcal{R}_{\rm upper}=  \demand^\ast-\frac{y^\ast}{a}$. Since 
$y^\ast=(1-\proportion)[1-2\Phi(\price)]$ then: 
\begin{equation}\label{lota}
\mathcal{R}_{\rm lower}\stackrel{\rm (\ref{tena})}{\geq} 
\frac{\price}{J} +\frac{1-\alpha}{\proportion \sus} \frac{a}{b}\price,\, \forall p\geq 0, \hbox{ and } 
\mathcal{R}_{\rm upper}\stackrel{\rm (\ref{estF})}{\leq} \frac{\price}{J} +\frac{1-\proportion}{\proportion\sus} \frac{b}{a} \price,\, \forall p\in [0, h].
\end{equation}

Next, due to (\ref{hhh}), if  $0\leq \price\leq h$ then 
$x\in [0,\, \frac{1-\proportion}{\proportion}\frac{b}{a} \price]$ implies that
$x\in [0,\, x_r]$. This implication and (\ref{vera}) ensure that
$2\proportion\sus \Phi^\prime(x)\geq a$, in case $x\in [0,\, \frac{1-\proportion}{\proportion}\frac{b}{a} \price]$ and $p\in [0, h]$. This conclusion and the relation (\ref{edive}) between $g_\price^\prime$ and $\Phi^\prime$ together ensure that if $p\in [0, h]$ then 
$g^\prime_\price(\demand)\geq a$ for each $\demand\in \left[\frac{\price}{\sus},\,\, 
\frac{\price}{\sus} + \frac{1-\proportion}{\proportion}\frac{b}{a} \price\right]$. 
From this inequality, it follows that the graph of $g_\price(\demand)$, 
$\demand\in \left[\frac{\price}{\sus},\, \frac{\price}{\sus}+ 
\frac{1-\proportion}{\proportion}\frac{b}{a} \price\right]$, lies above 
the line of slope $a$ constructed by us, provided $\price\in [0, h]$. This fact implies that
$\mathcal{R}(p)\leq \mathcal{R}_{\rm upper}$, when $\price\in [0, h]$ (for this implication 
to be valid it is important that the interval on which $g_\price(\cdot)$ lies above the line extends
up to the upper bound of $\mathcal{R}_{\rm upper}$ provided by (\ref{lota})). Then, via simple 
algebraic manipulation with the r.h.s. of the estimate (\ref{lota}) for 
$\mathcal{R}_{\rm upper}$, the assertion (b) of the lemma follows.

Finally, since $g_\price^\prime(\cdot)\leq b$ (as (\ref{gadiv}) ensures) 
then the graph of $g_\price(d)$, $d\in \big[\frac{\price}{\sus}, \infty\big)$, lies below 
the line with slope $b$ constructed by us. Hence, $\mathcal{R}_{\rm lower}\leq \mathcal{R}(\price)$, 
and the assertion (a) of the lemma follows  via simple algebraic manipulation with the 
r.h.s. of the lower bound for $\mathcal{R}_{\rm lower}$ provided by (\ref{lota}).\hfill$\Box$

\bigskip\noindent{\bf Lemma~\ammellabel{sitin}.}\ \ {\em
Let   $(\price_n,\demand_n), \  n\geq 0$,  be an orbit of  (\ref{DYMA}).
Then, either $(\price_n,\demand_n) \to (0,0)$ or there exists finite $r \geq 0$ 
such that  
\begin{equation}\label{ilu}
\hbox{either }0 \leq \price_{r+1} \leq \feedback \demand_r\hbox{ or } 
\feedback \demand_r\leq \price_{r+1} \leq   0.
\end{equation}
}

\medskip\noindent {\bf Proof.} \ {\em Case 1:} \ First, we consider the case when  
$\demand_0$ and  $\price_0$ have different signs. We shall conduct our arguments 
under the assumption that  $\price_0 \geq 0$ and  $\demand_0 \leq 0$; 
in the opposite case, i.e., when $\price_0 \leq 0$ and  $\demand_0 \geq 0$, the proof 
follows by exactly the same argument with the obvious change of inequality 
directions.

The assumption $\price_0 \geq 0$ and  $\demand_0 \leq 0$ and the relation (\ref{clan}) 
imply that: 
\begin{equation}\label{shal}
\price_1 \leq \price_0\hbox{ and }\price_1 \geq \feedback\demand_0.
\end{equation}
For this reason, if  $\price_1 \leq 0$ then (\ref{ilu})  holds with $r=0$ and thus, the lemma is proved.
Assume the contrary:  $\price_1 \geq 0$. Let us prove that $\demand_1\leq 0$ under this
assumption.  Since
\begin{equation}\label{qot}
|\demand_1 - \mathcal{D}(\price_1)|\stackrel{\rm(\ref{blan}), (\ref{kabe})}{=} 
|g_{\price_1}(\demand_0)-g_{\price_1}(\mathcal{D}(\price_1)|\stackrel{\rm(\ref{squez})}{\leq} 
|\demand_0 - \mathcal{D}(\price_1)|,
\end{equation}
then $|\demand_1 - \mathcal{D}(\price_1)|\leq |\demand_0 - \mathcal{D}(\price_1)|$. The 
inequality (\ref{qot}), the inequality $\mathcal{D}(\price_1)\leq 0$ (that stems from 
(\ref{pozd}) and the assumption $\price_1 \geq 0$) and the inequality $\demand_0\leq 0$ 
(the initial assumption)  
can hold simultaneously only if $\demand_1 \leq 0$. This is the conclusion we aimed for.

The argument of the above paragraph either finishes the proof of the lemma or 
implies  that $\price_1\geq 0$ and $\demand_1\leq 0$ hold. In the latter case,
the argument can be repeated with these  inequalities in the place of the 
argument's assumption ``$\price_0 \geq 0$ and  $\demand_0 \leq 0$''.  Upon 
repeating this procedure,  
we either eventually find an $r$ such that $\price_{r+1}\leq 0$ or not. When the first alternative
is the case, the last but one repetition of the argument ensured that $\price_{r}\geq 0$.
Combining this inequality with  $\price_{r+1}\leq 0$ and equation (\ref{clan}), we get 
that the second double inequality of (\ref{ilu}) holds true for the found value of $r$,
and thus the lemma's proof is finished. The second alternative is analyzed in 
the paragraph below. 

Let us assume that the argument of the last but one paragraph can be 
repeated infinitely. Since, after $n$-th repetition, we have that $\price_{n}\geq 0$ 
(otherwise, $r$ would have been equal to $n-1$) and $\demand_{n}\leq 0$ (this inequality is 
the final conclusion in each repetition) then, in virtue of (\ref{clan}), 
$\price_{n+1}\leq \price_n$ and therefore, the sequence $\{\price_n\}_{n=0}^\infty$ 
must have a non-negative limit; we denote it by $\bar{\price}$. The convergence  
$\price_n\rightarrow \bar{\price}$ and the relation 
$\demand_{n-1}=\feedback^{-1}(\price_n-\price_{n-1})$ (that follows from (\ref{clan})) 
ensure then that $\demand_n\rightarrow 0$. Therefore,  $(\bar{\price}, 0)$ is an 
equilibrium of (\ref{DYMA}). Due to Theorem~\ref{propdyns}, $\bar{\price}=0$. This proves 
the lemma's assertion.

\medskip\noindent{\em Case 2:} \ 
Now we consider the case where  $\demand_0$ and  $\price_0$ have the same sign.
We shall assume that  $\price_0 \geq 0$ and  $\demand_0 \geq 0$; the opposite 
case can be treated by exactly the same argument.

We observe initially that 
\begin{equation}\label{ore}
\hbox{if for some }n\geq 1,\,\,\, \demand_{n-1} \geq  0 \hbox{  and  }\price_{n-1}\geq 0
 \hbox{ then }
\price_{n}  \geq \price_{n-1}\hbox{ and }\demand_{n} \leq \demand_{n-1}. 
\end{equation}
The first conclusion in the implication (\ref{ore}) stems directly from its assumption 
and from (\ref{clan}). For the second one, we argue as follows. If 
$\mathcal{R}(\price_n)\not=\pm\infty$ then the following relations hold:
\begin{equation}\label{viol}
|\demand_n|\stackrel{\rm (\ref{blan}), (\ref{mabe})}{=}\left|g_{\price_n}(\demand_{n-1})-g_{\price_n}
\left(\mathcal{R}(\price_n)\right)\right|
\stackrel{\rm (\ref{squez})}{<}  \left|\demand_{n-1}-\mathcal{R}(\price_n)\right|.
\end{equation}
{}From (\ref{viol}), we obtain that $|\demand_n|<|\demand_{n-1}-\mathcal{R}(\price_n)|$. 
Next, from $\price_{n}\geq 0$ and (\ref{tvah}), we obtain that $\mathcal{R}(\price_n)\geq 0$. 
The two inequalities just obtained and the assumption $\demand_{n-1}\geq 0$ of (\ref{ore}) can hold altogether only if $\demand_{n}\leq \demand_{n-1}$. This completes the proof of (\ref{ore}) in case $\mathcal{R}(\price_n)\not=\pm\infty$. Let us analyse the opposite case. First, we note that
the inequality $\price_n\geq 0$ (which is the consequence of the assumption $\price_{n-1}\geq 0$ and our conclusion $\price_{n}\geq \price_{n-1}$) and the property (\ref{tvah}) altogether imply that 
$\mathcal{R}(\price_n)$ can only be $+\infty$. This conclusion and property (\ref{mabe}) 
ensure then that $g_{\price_n}(\demand)<0$ for all $\demand$, and consequently, that  
$\demand_{n}=g_{\price_{n}}(\demand_{n-1})<0$. Since $\demand_{n-1}\geq 0$ has been assumed then 
$\demand_n\leq \demand_{n-1}$ follows and (\ref{ore}) is established.  

Consider now an orbit $\{(\price_n, \demand_n)\}_{n=0}^\infty$ satisfying
our assumption $\price_0 \geq 0$ and  $\demand_0 \geq 0$. Obviously, either $\demand_n\geq 0$ 
for all $n$, or we shall find $m<\infty$ for which  $\demand_{m} \leq 0$; there may 
be several such numbers, let $m$ denote the smallest of them. For this $m$ it holds that 
$\price_m\geq 0$. This is ensured by recursive application of (\ref{ore}) for all $n\leq m$.
Hence, the portion $\{(\price_n, \demand_n)\}_{n=m}^\infty$ of the considered orbit 
can be treated by the argument of {\em Case 1}. By this argument, the lemma's assertion 
follows. It remains only to prove the lemma for the case when $\demand_n\geq 0$ for all $n$. 
We do this below.
 
Assume (in addition to the assumption $\price_0 \geq 0$ and  $\demand_0 \geq 0$ already made) 
that $\demand_n\geq 0$ for all $n$. The relation (\ref{ore}) implies that in this case, $\{\price_n\}$ 
is a non-decreasing sequence of non-negative numbers. If it is unbounded from above, then so is 
the sequence $\{\mathcal{R}(\price_{n})\}$ in virtue of Lemma~\ref{estrh}~(a). But since, 
by the very definition, $|\demand_n|\leq 1$ for all $n$, then there must exist $k$ such that $\demand_k<\mathcal{R}(\price_{k+1})$. For this $k$, 
\[
\demand_{k+1}\stackrel{\rm (\ref{blan})}{=}g_{\price_{k+1}}(\demand_k)\left\{\begin{array}{l}
\stackrel{{\rm (\ref{moveh}),}\,\,
\demand_k<\mathcal{R}(\price_{k+1})}{<}g_{\price_{k+1}}\left(\mathcal{R}(\price_{k+1})\right)\stackrel{\rm (\ref{mabe})}{=}0, \hbox{ when }\mathcal{R}(\price_{k+1})\hbox{ is finite},\\
<0,\hbox{ since }g_{\price_{k+1}}(\cdot)<0\hbox{ when }\mathcal{R}(\price_{k+1})=+\infty
\end{array}\right.
\] 
(we used above that $\mathcal{R}(\price_{k+1})$ cannot be $-\infty$ when $\price_{k+1}\geq 0$) 
and therefore, $\demand_{k+1}<0$. This contradicts the current assumption that
$\demand_n\geq 0$ for all $n$. Hence,  $\{\price_n\}$ must be bounded and thus, 
converges to a finite limit $\bar{\price}$. From this convergence we deduce that $(\price_n, \demand_n)\rightarrow (0,0)$ in the same way as it has been done four paragraphs above.  
\hfill$\Box$

\bigskip\noindent{\bf Lemma~\ammellabel{sitin1}.}\ \ {\em
Let\footnote{$b$, $a$ and $h$ employed here for the definition of $k$ and $\feedback_c$ 
                       have been defined in (\ref{debe}), (\ref{aaa}) and (\ref{hhh}).} 
$k$ be the minimal positive integer for which  $b^k \leq  \frac{1}{2}\frac{J\alpha (1-b)}{b}$ 
(the condition $b\in (0,1)$ guarantees the existence of such $k$).
Let $\feedback_c$ be defined as the maximal positive real number satisfying the following 
inequalities: 
\begin{equation}\label{cico}
(\fedum)\,\, \feedback_c k \leq  \frac{1}{2} \frac{ \sus\proportion (1-b)}{b},\,\,\,\, (\feddois)\,\,
\feedback_c\leq \proportion\sus \frac{a}{b},\,\,\,\, (\fedtres)\,\, \feedback_c\leq h.
\end{equation}

Suppose  $\feedback \in (0, \feedback_c]$. Let  $(\price_n, \demand_n), n\geq 0$, be an orbit of (\ref{DYMA}). Then,\newline 
if $0\leq \price_1\leq \feedback \demand_0$, then either $(\price_n, \demand_n)\rightarrow (0,0)$ or  
there exists a finite $\ell\geq 1$ such that 
\begin{equation}\label{api}\begin{array}{l}
\hbox{(i) }\demand_{j-1}  \geq \mathcal{R}(\price_{j}),\,\,\hbox{(ii) }0 \leq \demand_{j} \leq  b \demand_{j-1},\,\, \hbox{(iii) }\price_{j+1}  \geq \price_{j},\,
 \hbox{ for all }j =1,\cdots, \ell,\\
\hbox{ and (iv) }\demand_{\ell} \leq \mathcal{R}(\price_{\ell+1}); \end{array}
\end{equation}
if $\feedback \demand_0\leq\price_1\leq 0$, then either $(\price_n, \demand_n)\rightarrow (0,0)$ or  
there exists a finite  $\ell \geq 1$ such that 
\begin{equation}\label{rapi}\begin{array}{l}
\hbox{(i) }\demand_{j-1}  \leq \mathcal{R}(\price_{j}),\,\,\hbox{(ii) }b \demand_{j-1}\leq \demand_{j}
 \leq 0   ,\,\, \hbox{(iii) }\price_{j+1}  \leq \price_{j},\,
 \hbox{ for all }j =1,\cdots, \ell,\\
\hbox{ and (iv) }\demand_{\ell} \geq \mathcal{R}(\price_{\ell+1}). \end{array}
\end{equation}}

\medskip\noindent {\bf Proof.} \ We shall conduct the proof under the assumption 
$0\leq \price_1\leq \feedback \demand_0$; in the case $\feedback \demand_0\leq \price_1\leq 0$ 
the argument is analogous and will not be presented.

{\em Step 1:} \ We claim that the following inequality stems from lemma's assumptions:
\begin{equation}\label{has}
\demand_0\geq \mathcal{R}(\price_1).
\end{equation}
Indeed, by the very definition (\ref{aver}), $|\demand_0|\leq 1$ and hence
the assumption $0\leq \price_1\leq \feedback \demand_0$ implies that $0\leq \price_1\leq \feedback$. This, together with the assumption $\feedback\leq \feedback_c$ and the constraint (\ref{cico}-(\fedtres)) 
yield that 
$0\leq \price_1\leq h$ and thus we can apply Lemma~\ref{estrh}(b). It gives:
\begin{equation}\label{ovo}
\mathcal{R}(\price_1)\leq
\frac{1}{\proportion\sus}\frac{b}{a}\price_1.
\end{equation} 
On the other hand, the assumption 
$\price_1\leq \feedback\demand_0$  together with  the constraint (\ref{cico}-(\feddois)) imply that
$\demand_0\geq \frac{1}{\proportion\sus}\frac{b}{a}\price_1$. This inequality and (\ref{ovo}) yield (\ref{has}).

{\em Step 2:} \ We get:
\begin{equation}\label{lasa}
\demand_1\stackrel{\rm (\ref{blan})}{=}g_{\price_1}(\demand_0)\stackrel{\rm (\ref{moveh}),(\ref{has})}{\geq} g_{\price_1}\left(\mathcal{R}(\price_1)\right)\stackrel{\rm(\ref{mabe})}{=} 0,
\end{equation}
and 
\begin{equation}\label{mawa}
|\demand_1| \stackrel{\rm (\ref{blan}), (\ref{mabe})}{ = } \left|g_{\price_1}(\demand_0)-g_{\price_1}\left(\mathcal{R}(\price_1)\right)\right|
\stackrel{\rm (\ref{squez})}{\leq} b\left|\demand_0-\mathcal{R}(\price_1)\right| \leq b\demand_0, 
\end{equation}
where the last passage in (\ref{mawa}) is valid because $\demand_0\geq 0$ 
(this is one of the inequalities assumed in the beginning of the proof)  
and  $0\leq \mathcal{R}(\price_1)\leq \demand_0$ (here, the first inequality 
is provided by (\ref{tvah}) and the assumption that $\price_1\geq 0$, while 
the second inequality is identical to (\ref{has})); note also that this 
double inequality ensures that $\mathcal{R}(\price_1)\not=\pm\infty$ and hence, 
the last passage in both (\ref{lasa}) and (\ref{mawa}) is legitimate. 

{\em Step 3:} \ Since $\demand_1\geq 0$, as  ensured by (\ref{lasa}),  then (\ref{clan}) 
yields that  
\begin{equation}\label{basa}
\price_2\geq \price_1.
\end{equation}
This is the final conclusion of the third step.

The relations (\ref{lasa}), (\ref{mawa}) and (\ref{basa}) prove the validity of (i), (ii) and (iii) of (\ref{api}) for $j=1$. 

Now, if $\demand_1\leq \mathcal{R}(\price_2)$ then the lemma is proved with $\ell=2$. 
Thus, we continue the proof assuming the contrary:  $\demand_1> \mathcal{R}(\price_2)$. 
Obviously, this assumption implies that  $\demand_1\geq \mathcal{R}(\price_2)$. 
Taking this relation in the place of (\ref{has}) and repeating the second and the 
third steps of the argument presented above\footnote{Note that we do not need to repeat 
  the argument's first step since its conclusion 
  (\ref{has}) is now a direct consequence of our assumption. Hence, we will not 
  employ Lemma~\ref{estrh}(b) in the current and the consequent repetitions. 
  Getting rid of the necessity for the use of this lemma is here,
  because we cannot guarantee that $|\price_j|\leq h$, for $j\geq 2$, and therefore,
  we cannot ensure the validity of lemma's assumptions at $j$-th repetition for $j\geq 2$.}
we deduce the relations $0\leq \demand_2\leq b\demand_1$ and $\price_3\geq \price_2$. 
These relations and the inequality $\demand_1\geq \mathcal{R}(\price_2)$ 
imply that (i), (ii) and (iii) of (\ref{api}) hold for $j=2$. 

It is obvious that the argument of the above paragraph can be repeated for $j>2$ provided
$\demand_j>\mathcal{R}(\price_{j+1}$. 
However, the repetition process cannot last forever, unless $\price_1=0$. The reason 
for this is  the following. After $j$ consecutive repetitions, we would have that
$\demand_j\leq b^j\demand_0$ (in virtue of (\ref{api}-ii)) and $\price_{j+1}\geq 
\price_1\geq 0$ (in virtue of the assumption $\price_1\geq 0$ and the property 
(\ref{api}-iii)). Thus, the sequence $\{\demand_j\}$ decreases to zero, while the 
sequence $\{\mathcal{R}(\price_j)\}$ possesses -- in virtue of Lemma~\ref{estrh}(a) -- 
the following property: $\mathcal{R}(\price_{j+1})\geq \frac{1}{\proportion\sus}
\frac{a}{b}\price_{j+1}\geq \frac{1}{\proportion\sus}\frac{a}{b}\price_{1}$. Consequently, 
if $\price_1\not =0$ then there must be a finite $\ell$ for which $\demand_\ell\leq \mathcal{R}(\price_{\ell+1})$, and the lemma is proved.

To finish the proof, we have to complete the argument of the above paragraph by 
analyzing the case $\price_1=0$. Suppose in addition that $\demand_0=0$. It then 
follows directly from (\ref{DYMA}) that $(\price_j, \demand_j)=(0,0)\, \forall j\geq 1$. 
This conclusion completes the proof since one of the alternatives in the lemma's 
assertion is that $(\price_n, \demand_n)\rightarrow (0,0)$. 

The last case to be considered is thus, $\price_1=0$ and $\demand_0>0$. In this 
case, $g_{\price_1}(\cdot)\equiv g_0(\cdot)$. We note that (\ref{debe}) and 
(\ref{phicon}-\eqibda) ensure that $g_0(0){=}1-2\Phi(0){=}0$. Hence, we get:
\begin{equation}\label{danu}
\demand_1\stackrel{\rm(\ref{blan})}{=}
g_{\price_1}(\demand_0)\stackrel{\price_1=0}{=}g_0(\demand_0)\stackrel{\demand_0>0 {\rm\ and\ (\ref{moveh})}}{>}g_0(0)=0.
\end{equation}
{}From (\ref{clan})  and (\ref{danu}), we get that $\price_2=\feedback \demand_1>0$. 
This implies that $0\leq \price_2\leq \feedback\demand_1$. With the latter double 
inequality in the place of $0\leq \price_1\leq \feedback\demand_0$ we repeat the whole  
argument starting from the beginning of the proof and finishing at the end of the 
above paragraph. In the repetition, we shall not stumble upon the inconvenient 
possibility ``$\price_2=0$'', since  we have just shown that $\price_2>0$. The result is 
the following conclusion: either $(\price_n, \demand_n)\rightarrow (0,0)$ or
(\ref{api}) holds with $j=2, \ldots, \ell$ for some finite $\ell\geq 2$. As for
the validity of (\ref{api}-i, ii, iii) for $j=1$, it has been already established 
(read the sentence after (\ref{basa})). This completes the proof of the lemma. \hfill$\Box$
           
\bigskip\noindent{\bf Lemma~\ammellabel{sitin2}.}\ \ {\em
Suppose  $\feedback\in (0, \feedback_c]$ where $\feedback_c$ satisfies (\ref{cico}).  
Let $(\price_n, \demand_n), n\geq 0$, be an orbit of (\ref{DYMA}). Then,\newline
if $0 \leq \demand_0 \leq \mathcal{R}(\price_1)$ then 
either   $ (\price_n, \demand_n) \to  (0,0)$  or
there exists  $m \geq 1$ such that 
\begin{equation}\label{kawi}\begin{array}{l}
\hbox{(i) }\mathcal{D}(\price_1) \leq \demand_j \leq 0, 
\hbox{ for  } \ j=1, \cdots, m, \hbox{ (ii) }
0 \leq \price_j \leq \price_{j-1},\hbox{  for  }\ j=2, \cdots, m,\\
\hbox{ and (iii) }\feedback \demand_m \leq \price_{m+1 } \leq 0;\end{array}
\end{equation}
and if $0\geq \demand_0\geq \mathcal{R}(\price_1)$  then 
either   $ (\price_n, \demand_n) \to  (0,0)$ or 
there exists  $m \geq 1$ such that 
\begin{equation}\label{wawi}\begin{array}{l}
\hbox{(i) }\mathcal{D}(\price_1)\geq \demand_j\geq 0, 
\hbox{ for  } \ j=1, \cdots, m, \hbox{ (ii) }
0\geq \price_j \geq \price_{j-1},\hbox{  for  }\ j=2, \cdots, m,\\
\hbox{ and (iii) }\feedback \demand_m \geq \price_{m+1 } \geq 0.\end{array}
\end{equation}
} 

\medskip\noindent {\bf Proof.} \ We shall conduct the proof under the assumption 
$0 \leq \demand_0 \leq \mathcal{R}(\price_1)$; in the case 
$0\geq \demand_0\geq \mathcal{R}(\price_1)$ the argument is analogous 
and hence  will not be presented.

First, we conclude that  $\demand_1\leq 0$ reasoning as follows. If 
$\mathcal{R}(\price_1)\not=+\infty$, then
\begin{equation}\label{yola}
\demand_1\stackrel{\rm(\ref{blan})}{=}g_{\price_1}(\demand_0)
\stackrel{\rm (\ref{moveh}),\hbox{\small\ and the assumption }
\demand_0\leq \mathcal{R}(\price_1)}{\leq} g_{\price_1}
\left(\mathcal{R}(\price_1)\right)\stackrel{\rm (\ref{mabe})}{=}0.
\end{equation}
If $\mathcal{R}(\price_1)=+\infty$ then in virtue of (\ref{mabe}), $g_{\price_1}(\demand)<0$ for all $\demand$, and hence $\demand_1=g_{\price_1}(\demand_0)<0$.
 
Second, from (\ref{tvah}) and the assumption $\mathcal{R}(\price_1) \geq 0$ we get that 
$\price_1\geq 0$ and therefore, in virtue of Lemma~\ref{estAh}, $\mathcal{D}(\price_1)\leq 0$.   
The latter and the assumption $\demand_0\geq 0$ imply that $\demand_0\geq \mathcal{D}(\price_1)$.
We use this relation in the calculations below:
\begin{equation}\label{rato}
\demand_1\stackrel{\rm (\ref{blan})}{=} g_{\price_1}(\demand_0 )\stackrel{\rm (\ref{moveh}), 
\demand_0\geq \mathcal{D}(\price_1)}{\geq}  g_{\price_1}\left(\mathcal{D}(\price_1)\right)\stackrel{\rm(\ref{kabe})}{=}\mathcal{D}(\price_1),
\end{equation}
and we conclude thus, that $\demand_1\geq \mathcal{D}(\price_1)$.

At the third step, we turn our attention to $\price_2$, but consider separately the cases 
$\price_2\leq 0$ and $\price_2> 0$.  

Let us assume that $\price_2\leq 0$. Our aim is to prove that in this case  
(\ref{kawi}) is true for $m=1$. Since  it  has been proved in the previous two 
steps that  $\mathcal{D}(\price_1)\leq \demand_1\leq 0$  and since (\ref{kawi}-(ii)) 
is void for $m=1$, then our aim is achieved as soon as we prove that 
$\feedback \demand_1\leq \price_2\leq 0$. In this relation, the second inequality is 
exactly our current assumption, and the first inequality is valid because  
$\price_2=\price_1+\feedback \demand_1$ (the relation (\ref{clan})) and  because 
$\price_1\geq 0$ (proved at the second step above). The lemma is  therefore proved in 
case $\price_2\leq 0$. 

Let us now assume that $\price_2> 0$. Our objective is to show that (\ref{kawi}-(i)) 
and (\ref{kawi}-(ii))  are satisfied with $j=2$. 

We start by  establishing (\ref{kawi}-(ii)) for $j=2$. The  first 
inequality is a direct consequence of our current assumption. The  second inequality stems from
the relation $\price_2=\price_1+\feedback \demand_1$ (the eq.~(\ref{clan})) and 
the inequality $\demand_1\leq 0$ that has been proved in the first step.

We proceed by proving (\ref{kawi}-(i)) for $j=2$. Using $\demand_1\leq 0$ (proved in the first step) 
and $\price_2>0$ (the current assumption) we get:
\[
\demand_2\stackrel{\rm (\ref{blan})}{=}g_{\price_2}(\demand_1)\stackrel{{\rm (\ref{moveh}),}\,\,
\demand_1\leq 0}{\leq}g_{\price_2}(0)\stackrel{{\rm (\ref{moveh}), (\ref{tupo}),}\,\, \price_2>0}{\leq}g_{\price_2}\left(\frac{\price_2}{\sus}\right)\stackrel{\rm (\ref{tupo})}{\leq} 0,
\]
that proves the second inequality in  (\ref{kawi}-(i)) for $j=2$. 
Thus, to achieve our objective, we have only to show that 
\begin{equation}\label{ode}
\demand_2\geq \mathcal{D}(\price_1).
\end{equation}
We split the argument in two parts. First we prove (\ref{ode}) under the assumption 
that $\demand_1\geq \mathcal{D}(\price_2)$, and then under the assumption that 
$\demand_1\leq \mathcal{D}(\price_2)$. 

Assume that $\demand_1\geq \mathcal{D}(\price_2)$. This yields the desired relation (\ref{ode}) 
as follows: 
\begin{equation}\label{aqz}
\demand_2\stackrel{\rm (\ref{blan})}{=}g_{\price_2}(\demand_1)\stackrel{{\rm (\ref{moveh}),}\,\, 
\demand_1\geq \mathcal{D}(\price_2)}{\geq} g_{\price_2}\left(\mathcal{D}(\price_2)\right)\stackrel{\rm (\ref{kabe})}{=}\mathcal{D}(\price_2)
\stackrel{\rm (\ref{juti})\hbox{  and }0\leq \price_2\leq \price_1}{\geq}\mathcal{D}(\price_1)
\end{equation}
(note that the inequality $0\leq \price_2\leq \price_1$ used in the last passage is provided 
by (\ref{kawi}-(ii)) for $j=2$ that we have proved above).

Suppose now that $\demand_1\leq \mathcal{D}(\price_2)$. Then 
\begin{equation}\label{olha}
\demand_2\stackrel{\rm (\ref{blan})}{=}g_{\price_2}(\demand_1)
\stackrel{{\rm (\ref{moveh}),\hbox{\small\ and the assumption}}\,\, 
\demand_1\leq \mathcal{D}(\price_2)}{\leq} g_{\price_2}
\left(\mathcal{D}(\price_2)\right)\stackrel{\rm (\ref{kabe})}{=}\mathcal{D}(\price_2)
\end{equation}
that ensures that $\demand_2\leq \mathcal{D}(\price_2)$. Next, we note that the 
assumption $\price_2>0$ and the property (\ref{pozd}) imply that $\mathcal{D}(\price_2)\leq 0$. 
Because of this inequality, the following two relations are valid in the present case: 
$\demand_1\leq \mathcal{D}(\price_2)\leq 0$ and $\demand_2\leq \mathcal{D}(\price_2)\leq 0$. 
These relations and the inequality
\begin{equation}\label{acol}
\left| \demand_2-\mathcal{D}(\price_2)\right|\stackrel{\rm(\ref{blan})}{=}\left|g_{\price_2}(\demand_1)-g_{\price_2}\left(\mathcal{D}(\price_2)\right)\right|\stackrel{\rm (\ref{squez})}{<} \left|\demand_1-\mathcal{D}(\price_2)\right|
\end{equation}
can hold altogether only if $\demand_1\leq \demand_2$. The latter and the relation 
$\mathcal{D}(\price_1)\leq \demand_1$ imply (\ref{ode}).

We can now conclude that if $\price_2>0$ then (i) and (ii) of (\ref{kawi}) 
hold for $j=2$. Recall that in the first two steps of the proof we have already 
concluded that (\ref{kawi}-(i)) is valid for $j=1$. Combining these conclusions, 
we close our argument's third step: when $\price_2>0$ then 
(i) and (ii) of (\ref{kawi}) hold for $m=2$. 

 
The third step of our argument, i.e., the proof that started right after (\ref{rato}) 
and finished above, can be repeated for $\price_3, \price_4, $ etc. until 
we find $m$ such that $\price_{m+1}\leq 0$. The argument will then ensure that 
(\ref{kawi}-(iii)) holds for the value of $m$ founded, and the previous recursion steps 
will ensure that (\ref{kawi}-(i)) and (\ref{kawi}-(ii)) hold for this value of $m$. This
finishes the proof of the lemma in case $m$ is founded. If, to the contrary, such an 
$m$ does not exist then the recursion ensures -- via (\ref{kawi}-(ii)) -- that 
$\{\price_n\}_{n=1}^\infty$ is a monotone non-increasing sequence bounded below   
by $0$. Let $\bar\price$ denote its limit. Then, from the equation (\ref{clan})  we
get that $\demand_n=(\feedback)^{-1}\left(\price_{n}-\price_{n-1}\right)$, and therefore, $\demand_n\rightarrow 0$ as $n\rightarrow \infty$. But, in virtue of Theorem~\ref{propdyns}, 
the $n\rightarrow \infty$ limit of $(\price_n, \demand_n)$ can only be $(0,0)$. This 
completes the proof of the lemma.  \hfill$\Box$  
                  
\medskip
\noindent{\bf  Proof of Thm.~\ref{dyns}.} \ 
In virtue of Lemma~\ref{sitin}, in order to establish the theorem, it is sufficient to 
prove that $(\price_n, \demand_n)\rightarrow (0,0)$, in case when either 
$\feedback  \demand_0 \geq \price_1 \geq 0$ or $\feedback \demand_0 \leq \price_1 \leq 0$. 
We shall assume that
\begin{equation}\label{wali}
\feedback  \demand_0 \geq \price_1 \geq 0
\end{equation}
for the rest of the proof. In the  case $\feedback \demand_0 \leq \price_1 \leq 0$, 
the proof is the same up to obvious changes of the inequality directions.  

The assumption (\ref{wali}) allows us to apply  Lemma~\ref{sitin1}. 
Let  $\ell$ be the integer provided by it (if $\ell$ does not exist then 
$(\price_n, \demand_n)\rightarrow (0,0)$ in accordance to the lemma and thus, the 
theorem is proved). 

First, let us prove that
\begin{equation}\label{lot}
\price_\ell \leq \frac{1}{2} \frac{\sus \proportion  (1-b)}{b}\demand_0.
\end{equation}
For the proof, we shall need the fact that $\price_\ell\geq 0$. It follows from the 
assumption $\price_1\geq 0$ via a recurrent application of the relation (\ref{api}-iii) 
(that is valid due to Lemma~\ref{sitin1}). The continuation of our argument 
depends on whether $(\ell-1)\geq k$ or $(\ell-1)<k$, where, recall, $k$ has been 
defined in Lemma~\ref{sitin1} as the minimal positive integer for which $b^k\leq \frac{1}{2}\frac{\sus\proportion(1-b)}{b}$; note that $k$ depends solely on the 
parameters of the studied dynamical system, and hence can be applied to any its orbit. 

If $(\ell-1)\geq k$ we argue as follows. The inequality $\price_\ell\geq 0$ just proved and
the property (\ref{tvah}) yield the inequality $\price_\ell \leq \sus \mathcal{R}(\price_\ell)$, 
from which we get (\ref{lot}) via the following chain of estimates:
\[
\begin{array}{rcll}
\price_\ell \leq \sus \mathcal{R}(\price_\ell)&\leq & \sus\demand_{\ell-1}& \hbox{ (due to (\ref{api}-i) of Lemma~\ref{sitin1} for }j=\ell\hbox{)}\\
&\leq & \sus b^{\ell-1} \demand_0 &\hbox{ (from (\ref{api}-ii) of Lemma~\ref{sitin1} applied for }j=1, \cdots, \ell-1\hbox{)}\\
&\leq & \sus b^k \demand_0 &\hbox{ (because }(\ell-1)\geq k\hbox{ and }b<1\hbox{)}\\
&\leq & \frac{1}{2} \frac{\sus \proportion  (1-b)}{b}\demand_0 &\hbox{ (by the definition of }k\hbox{)}.
\end{array}
\]
If $(\ell-1) < k$, then  (\ref{lot}) follows from
 the following chain of inequalities:  
\[
\begin{array}{rcll}
\price_\ell &= &\price_1 + \feedback \demand_1 + \cdots +\feedback \demand_{\ell-1} & \hbox{ (by 
recurrent application of (\ref{clan}))}\\
&\leq  & \feedback\demand_0+ \feedback \demand_1 + \cdots +\feedback \demand_{\ell-1} & \hbox{ (using 
the first inequality of (\ref{wali}))}\\
&\leq & \ell  \feedback\demand_0 & \hbox{ (using (\ref{api}-ii) of Lemma~\ref{sitin1} for }
j=1,\cdots, \ell-1,\\
& & & \hbox{\ \ \ \ \ \ \ \ and taking into account that }b<1\hbox{)}\\
&\leq & k  \feedback\demand_0 & \hbox{ (because }(\ell-1)<k\hbox{ in the considered case)}\\
&\leq &\frac{1}{2} \frac{\sus \proportion (1-b)}{b}\demand_0& \hbox{ (by the condition (\ref{cico}-(\fedum)) on }\feedback\hbox{)}.
\end{array}
\]

Second, from the inequality $\price_\ell\geq 0$ just established and  Lemma~\ref{estAh} we get that 
\begin{equation}\label{seaf}
\mathcal{D}(\price_\ell) \geq  
- \frac{b}{\proportion (1-b)} \price_\ell.
\end{equation} 

Third, we combine the conclusions (\ref{lot}) and (\ref{seaf}) to  deduce that
\begin{equation}\label{roze}
\mathcal{D}(\price_\ell)\geq -\frac{1}{2}\demand_0.
\end{equation}

The fourth step  is based on the  double inequality $0\leq \demand_\ell \leq \mathcal{R}(\price_{\ell+1})$; 
its first part is ensured by (\ref{api}-ii) with $j=\ell$ and the second part by (\ref{api}-iv).
We take this inequality in the place of the assumption of Lemma~\ref{sitin2} and derive then from the 
lemma that either   $(\price_n, \demand_n) \to  (0,0)$ and therefore the theorem is proved, or  
\begin{eqnarray}
\hbox{there exists  }m \geq 1&\mbox{such that} & \mathcal{D}(\price_{\ell+1}) \leq \demand_j \leq 0,\hbox{ for  }
 \ j=\ell+1, \cdots, \ell + m, \label{waz}\\
&\hbox{and}& \feedback \demand_{\ell+m} \leq \price_{\ell+m+1 } \leq 0. \label{woz}
\end{eqnarray}  

Combining (\ref{roze}) and (\ref{waz}), we get that 
\begin{equation}\label{blo}
-\frac{1}{2}\demand_0\leq \demand_{\ell+j}\leq 0, \hbox{ for }j=1, \ldots, m.
\end{equation}

The inequalities (\ref{blo}) just derived, the inequalities (\ref{api}-ii) for $j=1, \ldots, \ell$, and 
the inequality $b<1$ lead altogether to the following conclusion:
\begin{equation}\label{aga}\begin{array}{l}
\hbox{in the block }\demand_1, \ldots, \demand_{\ell}, \demand_{\ell+1},\ldots, \demand_{\ell+m}
\hbox{ of the sequence }\{\demand_n\}_{n=0}^\infty,\\
\hbox{the values of the first }\ell (\ell\geq 1) \hbox{ members are between }0 \hbox{ and }\demand_0\\
\hbox{while the values of the last }m (m\geq 1) \hbox{ members are between }-\frac{1}{2}\demand_0 
\hbox{ and }0.\end{array}
\end{equation}

The argument that started at (\ref{wali}) and finished  at (\ref{aga}) 
can be repeated  with the inequality 
(\ref{woz}) in the place of (\ref{wali}). The respective conclusion is that the block 
$\demand_1,\ldots, \demand_{\ell+m}$ is followed
by another block -- of the size 2 at least -- that consists of two  non-empty parts such that: 
the members of the first part are all between $\demand_{\ell+m}$ and $0$, and the members of the
second part are all between $0$ and  $-\frac{1}{2}\demand_{\ell+m}$. This implies in  virtue
of the inequality $|\demand_{\ell+m}|\leq \frac{1}{2}\demand_{0}$ (that stems from (\ref{aga})), 
that the absolute value of each members of the block does not exceed $\frac{1}{2}\demand_{0}$.

It is obvious that the argument can be repeated yielding at the $i$-th 
step the conclusion that the members of the corresponding block of the sequence 
$d_1, d_2, \ldots$ do not exceed $\left(\frac{1}{2}\right)^{i-1}\demand_0$, in modulus. 
{}From this, $\demand_n\rightarrow 0$ follows. 
           
To complete the proof, it is only left to show that $\price_n \to 0$.
Suppose that this is not the case. Then, there exists $ \varepsilon > 0$ and 
a sequence of integers $\{n_j\}_{j\in \mathbb{N}}$ such that $n_j  \uparrow \infty$
and   $|\price_{n_j}|  \geq \varepsilon$ for each $j$. We  suppose, without loss of
generality that   $\price_{n_j}  \geq \varepsilon$.  Thus:
\begin{equation}\label{mita}
\demand_{n_j}\stackrel{\rm(\ref{blan})}{=} g_{\price_{n_j}}(\demand_{n_{j-1}})\stackrel{\rm (\ref{monip})}{\leq}
 g_{\varepsilon}(\demand_{n_{j-1}}). 
\end{equation}
The inequality (\ref{mita}), the continuity of $g_\price(\cdot)$,\footnote{$g_\price(\cdot)$ 
is continuous because it is differentiable -- see (\ref{edive}).} and our 
conclusion $\demand_n\rightarrow 0$ imply altogether that $0\leq g_\varepsilon(0)$. 
But, on the other hand,  (\ref{moveh}) and (\ref{tupo}) ensure that $0> g_\varepsilon(0)$.  
Hence, by contradiction, we conclude that $\price_n\rightarrow 0$.

We thus have proved that $\demand_n\rightarrow 0$ and $\price_n\rightarrow 0$ for arbitrary $\price_0$ and $\demand_0$. This is the theorem's assertion. \hfill $\Box$

\section{Application}\label{vmeste}

\subsection{Our HIA model for a single risky asset market}\label{model}
The bridge that leads from our results to their applications is grounded on 
the Heterogeneous Interacting Agent Model to be presented now. We shall  
use the abbreviation ``HIAM''. When we need to distinguish it from other HIAMs, we shall 
refer to it as ``our HIAM''. 

The model's {\em time set} is $\{0,1,2,\ldots\}$.
The model's parameter  denoted by $K$ and called {\em the number 
of  model's agents} is an arbitrary natural number. There are abstract 
{\em agents} in the model; they are labeled by the numbers $1,2,\ldots, K$. 
The other model's parameter, denoted by $\proportion$ and  called {\em the proportion 
of speculators}, is an arbitrary real number in the interval  $(0,1)$. It plays the 
following role in the model: the agents numbered $1,2,\ldots, \proportion K$ are called 
{\em speculators}, and the rest of the agent population, i.e., the the agents numbered 
$\proportion K+1, \proportion K+2,\ldots, K$,  are called {\em fundamentalists}.
The choice of these names will be justified a few paragraphs below.

To each time $n\geq 1$ and to each agent labeled $k$, the model associates a random 
variable denoted by $d_n(k)$. Its definition is as follows  (the quantities involved 
in (\ref{rules}) and (\ref{rulef}) will be specified below; in particular, note that 
$v_n(k)$ will be  a random variable -- this is the source of the randomness of $d_n(k)$):
\begin{equation}\label{rules}
	\decide_n(k) = \left\{\begin{array}{l}
		+1, \hbox{ if \ \ } \sus \bard_{n-1} - \big\{\price_n-(\medaval+\ve_n(k))\big\} > 0\\
		-1, \hbox{ if \ \ } \sus \bard_{n-1} - \big\{\price_n-(\medaval+\ve_n(k))\big\} \leq 0
	\end{array}\right.\kern1em\begin{array}{r}\hbox{if } k=1,2,\ldots, \proportion K,\\
\hbox{i.e., if $k$ is a speculator},\end{array}
\end{equation}
\begin{equation}\label{rulef}
	\decide_n(k) = \left\{\begin{array}{l}
		+1, \hbox{ if \ \ }  - \big\{\price_n-(\medaval+\ve_n(k))\big\} > 0\\
		-1, \hbox{ if \ \ }  - \big\{\price_n-(\medaval+\ve_n(k))\big\} \leq 0
	\end{array}\right.\kern1em\begin{array}{r}\hbox{if }k=\proportion K+1, 
	\proportion K+2,\ldots, K,\\
\hbox{i.e., if $k$ is a fundamentalist}.\end{array}\end{equation}
$\decide_n(k)$ is called {\em the decision of agent $k$ at time $n$}; when $\decide_n(k)$ 
is $+1$ we say that ``the agent wishes to buy an asset share", and when it is  $-1$ 
we say that ``he wishes to sell an asset  share".\footnote{In our HIAM, there is no 
        an equivalent of the market clearing condition. Therefore, the one who wishes 
        to buy/sell will not always accomplish this wish. By this reason, we 
        why we say ``wishes to buy/sell'' and not simply ``buys/sells''.} 
This ``buy/sell''\ interpretation helps to understand the link of the model with 
the asset market and the names that we shall attribute to the model's parameters 
and variables  introduced below.

Let us specify the quantities involved in (\ref{rules}) and (\ref{rulef}). 

First, 
\begin{equation}\label{aver}
\bard_n:=\frac{1}{K}\sum_{k=1}^K d_n(k), \kern1em n=1,2,\ldots,
\end{equation} 
and in accordance with this definition and with the interpretation of each 
$\demand_n(k)$, the real world analogue of $\bard_n$ is the population relative 
excess demand for the asset at time $n$; the name for the analogue will be the name for
$\bard_n$ in the model.  However,  when possible,  we shall shorten this name to just 
the {\em excess demand}. Observe that (\ref{aver}) defines $\bard_n$ for $n\geq 1$.
As for $\bard_0$, it is one of the model's parameters called {\em the initial excess demand}; 
its value range is $[-1,1]$. The value of $\bard_0$ as well as of $\price_0$ to be defined 
below, must be specified in order to ``switch on''  the model's evolution.

Second, $\price_n$ is the model's quantity that corresponds to the price of one share of
the asset in the market mimicked by the model.\footnote{Actually, it corresponds to the
          to the logarithm of the price of one share of asset. By using the logarithm, we 
          allow the model's counterpart of price be negative. This simplifies mathematical 
          treatment of the model. Nevertheless, the name for $\price_n$ in the model 
          will be ``price'' rather than ``logarithm of price''.}  
We call $\price_n$  {\em the asset price at time $n$}, or simply {\em price}. The value 
range of $p_n$ is $\mathbb{R}$ for each time $n=0,1,2,\ldots$.

The model's price evolves in time. The evolution rule is as follows:
$\price_0$ is one of the model's parameters (it is called {\em the price initial value}), 
and all others $\price_n$'s are determined by the  recurrent relation
\begin{equation}\label{prule}
\price_{n+1}=\price_n+\feedback \bard_n,\kern1em n=0,1,2,\ldots,
\end{equation}
where $\feedback$ is a positive real model's parameter called {\em the feedback of the 
excess demand on price increment}. The recurrent relation (\ref{prule}) called 
{\em price update rule} mirrors the demand-and-supply law of real markets. Indeed,
when $\bard_n$ is positive [resp., negative] -- which means that there are more buyers 
than sellers [resp., sellers than buyers] -- the rule rises [resp., lowers] the price. 

Third, we specify and interpret $\medaval$ and $\ve_n(k)$'s. $\medaval$ is a real value 
which is one of the model's parameters, and $\{v_n(k), \, n=1,2,\ldots,  \, k=1,2,\ldots, K\}$ 
are independent among themselves and of all other  random variables present in the
model,  with a  common distribution function that will be denoted by $\Phi$; 
this function is another model's parameter. We assume that $\Phi$ satisfies the
the conditions (\ref{phicon}).
The value $\medaval+\ve_n(k)$ is called {\em the evaluation of the asset's fundamental 
value by agent $k$ at time $n$}, or, simply, {\em individual evaluation}. The name reflects 
the quantity that we want to mimic by $\medaval+\ve_n(k)$ in our model, namely: what 
the agent $k$ thinks the asset worths at time $n$, or as in traditional economical 
terms, the asset's fundamental value evaluated by agent $k$ at time $n$. 

In the real world,  asset's fundamental value is usually calculated 
from the official balance reports of the company that issued the asset. However, 
an individual evaluation of this value may not coincide with the officially calculated 
one because an individual may have privileged information (about the future company 
development, say), or may use his own way to deduce the fundamental value from the 
report data. This situation is mirrored in our model in the following way: $\medaval$ 
corresponds to the official asset's fundamental value that may be simply the average 
taken over all individual 
evaluations, and  differences between this value and individual evaluations are modeled
by the random variables $\ve$'s with the distribution function $\Phi$. 
Accordingly,  the name for $\medaval$ in our model is {\em the market fundamental value 
of the asset}, and the name for $\Phi$ is {\em the distribution of deviation around $\medaval$ 
of individual evaluations of the fundamental value of the asset}.
  
Up to now, we have introduced and interpreted all the quantities involved in the 
decision rule (\ref{rulef}). This allows us to make the following observation: 
an agent that ``acts'' due to this rule, will wish to buy one asset share if he thinks 
that the current price (i.e., $\price_n$)  is less than what the asset share worths 
in accordance to his current evaluation (i.e., $\medaval+v_n(k)$); otherwise, he will 
wish to sell one asset share. The one who behaves in this way in the real market, is called 
{\em fundamentalist}. This explains why we use this name for our model agents that 
obey the rule (\ref{rulef}). 

Finally, we explain the role of $\sus$. This model's parameter is a positive real 
number that mirrors {\em the social susceptibility}, or simply {\em susceptibility} 
of market agents. It models ``susceptibility'' because its use in the rule 
(\ref{rules}) makes it to correspond to the strength with which the population 
excess demand influences an individual to align his decision in accordance  
with this excess. Now observe that the population excess demand is the average 
of decisions taken over the whole population and hence the alignment means the 
coincidence of an individual decision with that of the majority of the population. 
Hence, the ``susceptibility'' mirrored by $J$ has a ``social'' aspect of an 
individual behavior. 

The parameter $\sus$ is also interpreted as the {\em traders' speculative trend}. 
The justification for this is as follows. 
In our model, the population excess demand is proportional to the price trend 
(by the price trend we mean $\price_{n+1}-\price_{n}$, and the proportionality 
is yielded by (\ref{prule})). Consequently, the agent  who aligns his decision with the 
population excess demand acts as an asset market speculator (i.e., an 
individual that wishes to buy [resp., sell] 
when the trend indicates the price will rise [resp., fall]). This
fact explains why the model's agents that ``use'' the decision rule
(\ref{rules}) are called {\em speculators}. However, these agents are not pure speculators
because the rule (\ref{rules}) mixes the speculative behavior and the 
fundamentalist  behavior (explained two paragraphs above). In this mixture, the weight of 
the speculative component is given by $\sus$. This role of $\sus$ suggests its interpretation
as the market traders' speculative trend. 

\subsection{Why and how the studied dynamical system mimics
evolution of asset price and excess demand in asset markets}\label{dynas}
Corollary~\ref{whywhy} will answer the question posed  in this section title. It 
grounds upon Theorem~\ref{karov} and  Lemma~\ref{svaz}, both to be stated below. 
We could state that the convergence affirmed in the theorem holds true because
the Strong Law of Large Numbers ensures that $\bard_n$ converges to a constant, as 
the number of agents increases. However, it is impossible to construct a rigorous 
proof based solely on this law. Some adequate estimates for deviations are necessary.
The idea behind them is simple but the execution requires heavy notations and 
cumbersome calculations. An interested reader may find all this in \cite{prado-thesis}.
As for the proof of the lemma, it can be conducted easily by induction in time $n$, 
and hence will not be presented here.

\medskip\noindent{\bf Theorem~\mishlabel{karov}} \ (On convergence of trajectories of our HIAM
to orbits of a (\ref{DYMA})-like\newline\mbox{\kern8em} dynamical system).\newline
 \ {\em Choose a distribution function
$\Phi$  satisfying the conditions (\ref{phicon})  and real numbers
$\proportion\in (0,1)$, $\feedback>0$, $\medaval\in \mathbb{R}$, $\sus>0$, $\price_0\in \mathbb{R}$
and $\bard_0\in [-1,1]$. 

For each $K=1,2,\ldots$, and the entities chosen above construct the HIA model 
as defined in Section~\ref{model} and denote the  price and the excess demand at time $n$ of 
the constructed HIAM  by $\price^{st,K}_n$ and $\bard^{st,K}_n$, 
respectively.\footnote{The purpose of these   modifications is to emphasize  that the 
         price and the demand in the HIA model constructed in Section~\ref{model} are 
         \underline{stochastic} quantities and that they depend  on $K$, the model's 
         population size.}

Next, for $n\geq 1$, define the sequence of pairs of real numbers
$\left\{\left(\price_n\, ,\,\, \demand_n\right),\kern1em n=0,1,2\ldots\right\}$ 
so that $\demand_0=\bard_0$ and for $n\geq 1$, 
\begin{equation}\label{fir}\left\{\begin{array}{l}
\price_n=\price_{n-1}+\feedback \demand_{n-1}\\
\demand_n=\proportion\big[1-2\Phi\big(\price_n- \medaval- \sus \demand_{n-1}
\big)\big]+(1-\proportion)
\big[1-2\Phi\big(\price_n-\medaval\big)\big]\end{array}\right. 
\end{equation} 

Then, for any  $t\in \mathbb{N}$ and any $\varepsilon>0$, it holds that
\[
\lim_{K\rightarrow \infty}\pr\Big( \big|\bard^{st, K}_n-\demand_n\big|\leq \varepsilon,\hbox{ and }
\big|\price^{st, K}_n-\price_n\big|\leq\varepsilon, \,\,\hbox{ for all }n=0,1,2, \ldots, t\Big)=1.
\]
}

\medskip
It can be easily noted that (\ref{fir}), the dynamical system that figures in Thm.~\ref{karov}, 
is slightly different from (\ref{DYMA}), the dynamical system that we study in Section~\ref{results}. 
Nevertheless, as the following lemma reveals, trajectories of these systems have a very simple and 
rigid link. This link allows one to ``transport'' any result about (\ref{DYMA}) into the 
framework of (\ref{fir}).

\bigskip\noindent{\bf Lemma~\ammellabel{svaz}} \ (on the link between the dynamical 
systems (\ref{DYMA}) and (\ref{fir})).\newline {\em Suppose  the dynamical systems 
(\ref{fir}) and (\ref{DYMA}) have the same values for  the parameters $\proportion$, $\feedback$, $\sus$,
and the same distribution function $\Phi$  and let the value of the parameter $\medaval$
of (\ref{fir}) be arbitrary. Let us add the superscript $\ast$ to $\price_n$ and $\demand_n$  
in the system (\ref{fir}) so that they might be distinguished from $\price_n$ and $\demand_n$  
in  the system (\ref{DYMA}). Suppose  that $\price_0^\ast$ and $\demand_0^\ast$, 
the initial conditions for an orbit of (\ref{fir}), relate to $\price_0$ and $\demand_0$, 
the initial conditions for an orbit of (\ref{DYMA}), as follows: \ \ 
$\price_0=\price_0^\ast-\medaval, \hbox{\ and\ \ }\demand_0=\demand^\ast_0$.

Then, for every $n\geq 1$, it holds that \ $\price_n=\price_n^\ast-\medaval \hbox{\ \ and\ \ }\demand_n=\demand^\ast_n$.
}

\bigskip
Thm.~\ref{karov} and Lemma~\ref{svaz} yield straightforwardly
 the following result. 

\bigskip\noindent{\bf Corollary~\corlabel{whywhy}} \ (why and how the 
dynamical system (\ref{DYMA}) mimics evolution of asset price and excess demand).\newline 
{\em Let $\price_0\in \mathbb{R}$, $\demand_0\in [-1,1]$, $\proportion\in (0,1)$, $\sus>0$, 
$\feedback>0$, $\medaval\in \mathbb{R}$, and a probability distribution function $\Phi$ 
satisfying the conditions (\ref{phicon}) correspond to characteristics of a single risky 
asset market as specified in the Introduction. 

Consider the dynamical system (\ref{DYMA}) determined by  $\proportion$, $\sus$, 
$\feedback$ and $\Phi$. For arbitrarely fixed time $T$, let $\left\{(\price_n, \demand_n)\, 
n=0,1,\ldots, T\right\}$ denote the orbit of the system starting from $(\price_0,\demand_0)$.

Then, the trajectory on the time interval $[0,T]$ of the difference between the price of 
the market asset and the asset's fundamental value lies closely to $\left\{\price_n, \, 
n=0,1,\ldots, T\right\}$, while the trajectory on the time interval $[0,T]$ of the 
excess demand for the asset lies closely to $\left\{\demand_n,\, n=0,1,\ldots, T\right\}$.  
Moreover, the larger the market agent population size the closer are the trajectories 
to their respective approximations.}

\subsection{Revealing and explaining 
market price and population excess demand dynamics properties}\label{application}
In the present section, we employ Corollary~\ref{whywhy} 
to interpret the results of Section~\ref{results} in terms of a real world asset market.  
For the interpretations to be valid, it is necessary that $T$ and the market agent population size  
(both from Corollary~\ref{whywhy}) be large. We take it for granted that both are 
as large as necessary, when necessary.

Before we proceed, it must be noted that assumptions of the results of 
Section~\ref{results} do not concern the whole distribution function 
$\Phi$ but rather solely  $\Phi^\prime(0)$ (that is, the value of 
the derivative of $\Phi$ at $0$). This is an interesting virtue of our results. However, 
when one wants to apply them to explain real world market properties this virtue 
becomes an  obstacle since there is no a generic interpretation for $\Phi^\prime(0)$ 
in terms of real world markets (that is, an interpretation that would suit 
any $\Phi$). To overcome this obstacle, we shall accept the following assumption:
\begin{equation}\label{predpo}
\Phi \hbox{ is a Normal distribution function with zero mean}. 
\end{equation}
The point here is that if $\sigma^2$ denotes the variance of $\Phi$ then
this assumption implies that $\sigma^2=\left(\sqrt{2\pi}\Phi^\prime(0)\right)^{-1}$ which, in turn,
allows us for the following interpretation:
\begin{equation}\label{perush}\begin{array}{l}
\left(\Phi^\prime(0)\right)^{-1}\hbox{ corresponds to the heterogeneity
of the distribution of individual}\\
\hbox{evaluations of the asset's fundamental value, in the sense that larger value of}\\ 
\hbox{$\left(\Phi^\prime(0)\right)^{-1}$ corresponds to higher dispersion of evaluations over the population}.\end{array}
\end{equation}
The interpretation (\ref{perush}) is extremely convenient in re-phrasing the 
results of Section~\ref{results} in terms of real world asset markets. This is a 
strong motivation for accepting the assumption (\ref{predpo}). There are other 
motivations. One of them stems from our belief that an individual is influenced by 
diverse factors and information streams, when he makes up his mind in respect 
to the fundamental value of an asset. Consequently, due to the Central Limit 
Theorem, individual evaluations of this value should be distributed over a population 
in accordance to a Normal Law. Another motivation comes from the fact that the Normal Law
would be the most convenient and robust theoretical model for a populational 
distribution obtained by sampling from a real world population. In other words, the 
Normal Law would likely be chosen for $\Phi$ when fitting our HIAM to a real world 
socio-economic process.

\bigskip\noindent{\bf Corollary~\corlabel{ahad}} (of Thm.~\ref{propdyns} obtained with the help of Cor.~\ref{whywhy}).\newline
{\bf (a)} \  {\em
The market state in which 
\begin{equation}\label{estate}\begin{array}{l}
\hbox{the market asset price coincides with the asset's fundamental value}\\
\hbox{and the excess demand for the asset is zero}
\end{array}\end{equation} 
is the only possible equilibrium state for the asset price and the excess demand 
of an asset market. 

This means, in particular, that if the asset price and the excess demand converge 
to some value as time goes on, then their limit values are respectively, 
the asset's fundamental value and $0$.}

\noindent{\bf (b)} (valid under additional assumption (\ref{predpo})) \ {\em 
Whether the state (\ref{estate}) is a locally stable or an unstable equilibrium is determined 
by a relation between two expressions that we present below and denote by $u$ and $w$:
\begin{equation}\label{pare}
\begin{array}{rcl}
u&:=&2\times\hbox{(proportion of speculators)}\times\hbox{(traders' speculative trend)}\times\\
&&\mbox{\kern2em}\times\hbox{(the heterogeneity of
the individual evaluations}\\
&&\mbox{\kern9em}\hbox{ of the asset's fundamental value)}{}^{-1}\\
\hbox{and}&&\\ 
w&:=&2\times\hbox{(the feedback of the excess demand on price increment)}\times\\
&&\mbox{\kern2em}\times\hbox{(the heterogeneity of
the individual evaluations}\\
&&\mbox{\kern9em}\hbox{ of the asset's fundamental value)}{}^{-1}.
\end{array}
\end{equation}
The determinant relation is as follows (it is illustrated in Figure~\ref{regions}): if 
\begin{eqnarray*}
w\leq 2u+2&& \hbox{and }u\leq 1\\ 
&&\hbox{(together with $u> 0, w> 0$ which are implicit from (\ref{pare}) and the natural }\nonumber\\
&&\hbox{constraints that impose that each entity in the definition (\ref{pare}) is positive)}\nonumber
\end{eqnarray*}
then (\ref{estate}) is locally stable, otherwise, it is unstable.
}

\bigskip\noindent{\bf Corollary~\corlabel{shtaim}} (of Thm.~\ref{dyns}  obtained with the help of Cor.~\ref{whywhy})).\newline 
{\em  
There exists a strictly positive threshold $\feedback_0$ such that 
when the feedback of the excess demand on the price increment 
is smaller than $\lambda_0$ then the local equilibrium described in Corollary~\ref{ahad} 
turns to be global, that means that
the long time limit values for the asset price and the excess demand 
for this asset are, respectively, the asset's fundamental value and $0$, whatever the 
initial values are.
}

\bigskip\noindent{\bf Corollary~\corlabel{arba}} (of Thm.~\ref{limo} and Conjecture~\ref{checked}  obtained with the help of Cor.~\ref{whywhy} and valid under additional assumption (\ref{predpo}).\newline
{\em 
Suppose that, as a result of a smooth change of the parameters' values, 
an asset market is transferred from the parameter values region in which its equilibrium is stable 
to the region in which its equilibrium is unstable. Suppose that the change is such that 
\begin{description}
\item{(i)}{\vskip-1cm} 
\begin{description}{\vskip-1cm}
\item{(a)}  the traders' speculative trend (i.e., the value of parameter $\sus$) increases, or 
\item{(b)}   the proportion of speculators (i.e., the value of parameter $\proportion$) increases, or 
\item{(c)}   the heterogeneity of the distribution of the individual evaluations of the asset's fundamental value  diminishes (i.e., the value of $\Phi^\prime(0)$ increases), or
\item{(d)}  any combination of (a)--(c).  
\end{description}
\item{(ii)} the feedback of the excess demand on the price increment (i.e., the parameter $\feedback$) 
is kept fixed, and its value is sufficiently small so that $\feedback/(\proportion \sus)< 2$ 
for all values of $\proportion$ and $\sus$ during the  process of change.
\end{description}

Then, in the transfer course, there will be a (maybe short) interval of time during which the market asset price and excess demand exhibit regular oscillation with non damped amplitude.
}

\bigskip We would like to close the presentation with several comments:\newline
{\bf (a)} \  The last statement of Corollary~\ref{ahad} is important 
because a real world asset market always experiences endogenous and/or exogenous
shocks, and therefore the  stability properties  of the equilibrium (\ref{estate}) become 
an important issue.\newline 
{\bf (b)} \ In real world markets the price update increments are small, hence 
Corollary~\ref{shtaim} suggests that, in the real world, when the state 
(\ref{estate}) is stable, it is actually globally stable.\newline
{\bf (c)} \ One can conclude from Theorem~\ref{limo} and Conjecture~\ref{checked}  
that for a specific combination of market parameters' values, there will appear 
oscillations of asset price and excess demand for that asset. This conclusion
has limited application since the precise estimation of market's parameters is 
a very difficult task. Contrasting, the conclusion  presented in Corollary~\ref{arba} 
has a practical value because it describes a qualitative aspect of market 
behavior. {}From this description, one indeed can derive useful results, like 
the following one: Let one know that the corollary's assumption (i), (ii) hold true and 
 suppose one observes that a market exhibits oscillations of the asset price and the excess demand.
Then, even though the oscillations might have disappeared, one can affirm that the market 
is in an unstable state.

\end{document}